\pdfoutput=1
\documentclass[thmsa,notitlepage,leqno]{article}%
\usepackage{amssymb}
\usepackage{amsfonts}
\usepackage[leqno]{amsmath}
\usepackage{sw20jart}
\usepackage{graphicx}
\usepackage[letterpaper,margin=1in]{geometry}%
\providecommand{\U}[1]{\protect\rule{.1in}{.1in}}
\begin{document}

\author{Pablo Azcue and Nora Muler\\Departamento de Matematicas, Universidad Torcuato Di Tella.}
\title{A multidimensional problem of optimal dividends with irreversible switching: a
convergent numerical scheme}
\date{}
\maketitle

\begin{abstract}
In this paper we study the problem of optimal dividend payment strategy which
maximizes the expected discounted sum of dividends to a multidimensional set
up of $n$ associated insurance companies where the surplus process follows an
$n$-dimensional compound Poisson process. The general manager of the companies
has the possibility at any time to exercise an irreversible switch into
another regime; we also take into account an expected discounted value at
ruin. This multidimensional dividend problem is a mixed singular
control/optimal problem. We prove that the optimal value function is a
viscosity solution of the associated HJB equation and that it can be
characterized as the smallest viscosity supersolution. The main contribution
of the paper is to provide a numerical method to approximate (locally
uniformly) the optimal value function by an increasing sequence of sub-optimal
value functions of admissible strategies defined in an $n$-dimensional grid.
As a numerical example, we present the optimal time of merger for two
insurance companies.


\end{abstract}

\section{Introduction}

In this paper we study the problem of the optimal dividend payment strategy
which maximizes the expected discounted sum of dividends, in a
multidimensional setup of $n$ associated insurance companies. We assume that
the surplus process follows a multidimensional compound Poisson process. The
general manager of the companies has the possibility to exercise an
irreversible switch into another regime at any time; we also take into account
an expected discounted value at ruin, which is due at the first time of ruin
of one of the companies, and may depend on the value of the surplus of all the
companies both at and before this time of ruin. This ruin value is a
generalization to the multidimensional setting of the Gerber-Shiu penalty
functions introduced in Gerber and Shiu \cite{GS 1998}.

The problem of optimal dividend payments in the one-dimensional case was
proposed by de Finetti \cite{De Finetti} and it was studied in different model
setups. In the compound Poisson risk model, this problem was studied by Gerber
\cite{Gerber} using a limit of an associated discrete problem, and by Azcue
and Muler \cite{AM 2005} using a dynamic programming approach; see also an
overview on this problem in Schmidli \cite{Schmidli book 2008} and in Azcue
and Muler \cite{AM Libro}. For the limit diffusion approximations, see for
example Asmussen and Taksar \cite{Asmussen Taksar 1997} and for spectrally
negative L\'{e}vy risk processes see, for instance, Avram, Palmowski and
Pistorious \cite{APP 2007} and Loeffen \cite{Loeffen2008}.

In the one dimensional case, the final value of the portfolio at ruin is
non-positive and it is called a penalty. Let us mention for instance Dickson
and Waters \cite{DicksonWaters2004}, where the shareholders take care of the
deficit at ruin; Gerber, Lin and Yang \cite{GerberLinYang2006} where the
penalty is a function depending on the deficit at ruin; Thonhauser and
Albrecher \cite{TA} where they address the optimal dividend problem with
constant penalty. The optimal dividend problem in the spectrally negative
L\'{e}vy setting was solved by Loeffen and Renaud \cite{LoeffenRenaud2010}
with an affine penalty function, and by Avram, Palmowski and Pistorius
\cite{APP 2015} with a general penalty function depending on the deficit at ruin.

The one dimensional dividend problem with the possibility of an irreversible
switch was addressed by Ly Vath, Pham and Villeneuve \cite{LPV} in the
Brownian motion setup and by Azcue and Muler \cite{AM Switching} in the
compound Poisson setting.

The problem of dividend payment in the case of two insurances companies in the
compound Poisson risk model was studied by Czarna and Palmowski \cite{CP} for
a particular dividend strategy of reflecting two-dimensional risk process from
the line, and by Albrecher, Azcue and Muler \cite{AlAZMU} where they study the
optimal dividend strategy for two collaborating companies.

In this paper, the multidimensional dividend problem is a mixed singular
control/optimal problem. Its associated Hamilton-Jacobi-Bellman equation (HJB)
involves a first-order integro-differential operator, an obstacle operator and
$n$ derivative constraints; the integro-differential operator corresponds to
the discounted infinitesimal operator of the compound Poisson process, the
obstacle is related to the value of the portfolio after switching and the
derivative constraints are related to the dividend payments of the companies.
We prove that the optimal value function is a viscosity solution of the HJB
equation, that it can be characterized as the smallest viscosity supersolution
and also that a convergent limit of a family of admissible strategies that is
a viscosity solution of the associated HJB equation should be the optimal
value function (verification result). These results are natural extensions of
the results of \cite{AM Switching} to the multidimensional setting.

The way in which the optimal value function\ solves the HJB equation in the
$n$-dimensional state space suggests the optimal local control: in the closed
set where the optimal value function coincides with the obstacle
(\textit{switch region}), an immediate switch should be done; in the interior
of the set where the integro-differential operator is zero
(\textit{non-action} \textit{region}), no dividends are paid; and in the
interior of the set in which one or more of the derivative constraints are
tight (\textit{dividend payment region}), the corresponding companies pay a
lump sum of dividends. However, it is not clear what the optimal local control
is in the \textit{free boundaries} between the non-action region and the
dividend payment region. In the one dimensional case the "free boundaries" are
indeed "free points", and it can be seen that the optimal local control at
these points is just to pay all the incoming premium as dividends, so the
control surplus stays there until the arrival of the next claim. This is the
reason why the optimal strategy has a band structure and this free points can
be obtained by one-dimensional optimization techniques, see \cite{AM
Switching}. It is a hard task to obtain the free boundaries in the
multidimensional setting and there is no hope of finding a closed-form
solution for the optimal value function. The main contribution of this paper
is to provide a numerical method to approximate (locally uniformly) the
optimal value function by a sequence of sub-optimal value functions of
admissible strategies defined in an $n$-dimensional grid. These sub-optimal
value functions solve a discrete version of the HJB equation, and the
corresponding sub-optimal strategies are constructed partitioning the grid in
switch, non-action and dividend payment regions; so we also obtain numerical
approximations of the optimal switch, non-action and dividend payment regions
and the free boundaries between them.

For a convergence analysis of a numerical scheme for multidimensional singular
control problems in the diffusion setting using Markov chain approximation
methods, let us mention Kushner and Martins \cite{KM} and Budhiraja and Ross
\cite{BR}; see also the book of Kushner and Dupuis \cite{KD} for an exhaustive
survey. Regarding convergence of numerical schemes using the viscosity
solution approach, let us mention for instance Souganidis \cite{S} and Barles
and Souganidis \cite{BS}, where they propose a numerical scheme for
non-singular control problems in the context of the diffusion setting; roughly
speaking, they prove that the solutions of the numerical scheme converge to a
viscosity solution of the associated HJB equation and then, using a uniqueness
argument, they obtain the convergence result. In the numerical method of the
present work, there is not uniqueness of viscosity solutions in the HJB
equation; nevertheless, we construct numerically an increasing sequence of
value functions of a family of admissible strategies whose limit is a
viscosity solution of the associated HJB equation; then, using the
verification result mentioned above, we deduce that this limit is the optimal
value function.

As an application, we present the optimal time of merger (as change of regime
at the switch time) for two insurance companies. We show examples where the
non-action region could be non-connected even for exponential claim size
distributions. For a criteria of merger being an advantage over keeping the
two stand-alone companies under barrier strategies see Gerber and Shiu
\cite{GS Merger}.

The rest of the paper is organized as follows. In Section 2, we introduce the
model and derive some basic properties of the optimal value function. In
Section 3, we show that the optimal value function is a viscosity solution of
the corresponding (HJB) equation; we also characterize it as the smallest
viscosity supersolution and give a verification result. In Section 4, we
construct a family of admissible strategies at any point in a suitable grid.
In Section 5, we show that the discrete scheme convergences locally uniformly
by taking a suitable sequence of embedded grids. In Section 6, we present
examples of the problem of optimal merger time. Finally, in Section 7, there
is an Appendix with the proofs of the technical lemmas.

We use the following notation: $\mathbf{R}_{+}^{n}=[0,\infty)^{n}$, $\leq$
refers to the element-wise order on $\mathbf{R}^{n}$, $\mathbf{1}%
=(1,1,\ldots,1)\in$ $\mathbf{N}^{n}$, $\left(  \mathbf{e}_{i}\right)
_{i=1,...,n}\ $is the standard basis of $\mathbf{R}^{n}$, $\left[
\mathbf{x},\mathbf{y}\right]  =\left\{  \mathbf{z}\in\mathbf{R}^{n}%
:\mathbf{x}\leq\mathbf{z}\leq\mathbf{y}\right\}  $, $\mathbf{x}\vee
\mathbf{y}=(x_{1}\vee y_{1},...,x_{n}\vee y_{n})$, $\mathbf{x}\wedge
\mathbf{y}=(x_{1}\wedge y_{1},...,x_{n}\wedge y_{n})$.

\section{Model \label{Seccion Modelistica}}

Let us consider that the surplus process of $n$ companies, or branches of the
same company, follows an $n$-dimensional compound Poisson process with drift,
that means that the uncontrolled process $\mathbf{X}_{t}\in$ $\mathbf{R}%
_{+}^{n}$ can be written as
\begin{equation}
\mathbf{X}_{t}=\mathbf{x}^{0}+\mathbf{p}t-%
{\displaystyle\sum\nolimits_{k=1}^{N_{t}}}
\mathbf{U}_{k}.\label{UncontrolledSurplusOriginal}%
\end{equation}
Here $\mathbf{x}^{0}\in$ $\mathbf{R}_{+}^{n}$ is the initial surplus,
$\mathbf{p}=(p_{1},...,p_{n})$ where $p_{i}>0$ is the premium rate of company
$i$, $N_{t}$ is a Poisson process with intensity $\lambda$ and the downward
jumps $\mathbf{U}_{k}\in\mathbf{R}_{+}^{n}$ are i.i.d. vector random vectors
with joint multivariate distribution function $F$. We also assume that
$\mathbb{E(}\left\Vert \mathbf{U}_{k}\right\Vert \mathbb{)<\infty}$ and
$F(\mathbf{0})=0$. We call $\tau_{k}$ the time of arrival of the $k$-th jump
of the process, so $N_{t}=\max\{k:\tau_{k}\leq t\}$.

We can describe this model in a rigorous way by defining its filtered
probability space $(\Omega,\mathcal{F},\left(  \mathcal{F}_{t}\right)
_{t\geq0},\mathbb{P})$, where
\[
\Omega=\{(\tau_{k},\mathbf{U}_{k})_{k\geq1}\in\left(  \lbrack0,\infty
)\times\mathbf{R}_{+}^{n}\right)  ^{\mathbf{N}}:\tau_{k}<\tau_{k+1}\}
\]
and $\mathcal{F}_{t}$ is the $\sigma$-field generated by the set $\{(\tau
_{k},\mathbf{U}_{k}):\tau_{k}\leq t\}$. The uncontrolled surplus process
$\mathbf{X}_{t}$ is an $\mathcal{F}_{t}$-adapted c\`{a}dl\`{a}g (right
continuous with left limits) stochastic process. Each company pays dividends
to the same shareholders, let $\mathbf{L}_{t}\in$ $\mathbf{R}_{+}^{n}$ be the
vector of cumulative amount of dividends paid out up to time $t$ by each
company; we say that the dividend payment strategy $\mathbf{L}_{t}$ is
admissible if it is a non decreasing process, c\`{a}dl\`{a}g, adapted with
respect to the filtration $\left(  \mathcal{F}_{t}\right)  _{t\geq0}$ and
satisfies $\mathbf{L}_{0}\geq0$ and $\mathbf{L}_{t}\leq\mathbf{X}_{t}$ for any
$0\leq t<\tau^{\mathbf{L}}$, where $\tau^{\mathbf{L}}$ is the time in which
the process exits the set $\mathbf{R}_{+}^{n}$ due to a jump, that is
\begin{equation}
\tau^{\mathbf{L}}:=\inf\{\tau_{k}:\mathbf{X}_{\tau_{k}}-\mathbf{L}_{\tau
_{k}^{-}}\notin\mathbf{R}_{+}^{n}\}\text{.}\label{Definicion Tau L}%
\end{equation}
We define the \textit{controlled} surplus process as
\begin{equation}
\mathbf{X}_{t}^{\mathbf{L}}:=\mathbf{X}_{t}-\mathbf{L}_{t}\text{.}\label{XL}%
\end{equation}
It is not possible to pay any dividends once the controlled process
$\mathbf{X}_{t}^{\mathbf{L}}$ exits $\mathbf{R}_{+}^{n}$ so we extend
$\mathbf{L}_{t}=\mathbf{L}_{\tau^{\mathbf{L}}{}^{-}}$ for $t\geq
\tau^{\mathbf{L}}$. Note that $\mathbf{X}_{\tau^{\mathbf{L}}}^{\mathbf{L}%
}=\mathbf{X}_{\tau^{\mathbf{L}}{}^{-}}^{\mathbf{L}}-\mathbf{U}_{k_{0}}$ if
$\tau^{\mathbf{L}}=\tau_{k_{0}}$. At time $\tau^{\mathbf{L}}$, the
shareholders pay a penalty $\upsilon(\mathbf{X}_{\tau^{\mathbf{L}}{}^{-}%
}^{\mathbf{L}},\mathbf{U}_{k_{0}})$ (or get a reward in the case that
$\upsilon(\mathbf{X}_{\tau^{\mathbf{L}}{}^{-}}^{\mathbf{L}},\mathbf{U}_{k_{0}%
})$ is negative) depending on the surplus prior to ruin $\mathbf{X}%
_{\tau^{\mathbf{L}}{}^{-}}^{\mathbf{L}}$ and the size $\mathbf{U}_{k_{0}}$ of
the last jump of the uncontrolled process. Denote
\begin{equation}
B=\{(\mathbf{x},\mathbf{\alpha})\in\mathbf{R}_{+}^{n}\times\mathbf{R}_{+}%
^{n}\text{ s.t. }\mathbf{x}-\mathbf{\alpha}\notin\mathbf{R}_{+}^{n}%
\},\label{Definicion B}%
\end{equation}
the function $\upsilon:B\rightarrow\mathbf{R}$ generalizes the concept of
penalty at ruin. It is natural to assume that the penalty function
$\upsilon(\mathbf{x},\mathbf{\alpha})$ is non-increasing on $\mathbf{x}$ and
non-decreasing on $\mathbf{\alpha}$; furthermore, we assume that
$\mathbb{E}\left(  \left\vert \upsilon(\mathbf{0},\mathbf{U}_{1})\right\vert
\right)  <\infty$. The manager of the company also has the possibility at any
time $0\leq t<\tau^{\mathbf{L}}$ to exercise an irreversible switch whose
value is associated to a given function $f:\mathbf{R}_{+}^{n}\rightarrow
\mathbf{R} $. We assume that the function $f$ is either right continuous and
non decreasing or continuous.

Given an initial surplus $\mathbf{x}\geq0$, let us denote by $\Pi_{\mathbf{x}%
}$ the set of all pairs $\pi=(\mathbf{L},\overline{\tau})$ where $\mathbf{L}$
is an admissible dividend payment strategy and $\overline{\tau}$ is a switch
time\textit{. }We define%
\begin{equation}%
\begin{array}
[c]{lll}%
V_{\pi}(\mathbf{x}) & = & \mathbb{E}_{\mathbf{x}}\left(  \int_{0^{-}%
}^{\overline{\tau}\wedge\tau^{\mathbf{L}}}e^{-cs}\mathbf{a}\cdot
d\mathbf{L}_{s}+I_{\{\overline{\tau}<\tau^{\mathbf{L}}\}}e^{-c\overline{\tau}%
}f(\mathbf{X}_{\overline{\tau}}^{\mathbf{L}})\right) \\
&  & -\mathbb{E}_{\mathbf{x}}\left(  I_{\{\overline{\tau}\geq\tau^{\mathbf{L}%
}\}}e^{-c\tau^{\mathbf{L}}}\upsilon(\mathbf{X}_{\tau^{\mathbf{L}}{}^{-}%
}^{\mathbf{L}}{},\mathbf{X}_{\tau^{\mathbf{L}}{}^{-}}^{\mathbf{L}}%
{}-\mathbf{X}_{\tau^{\mathbf{L}}}^{\mathbf{L}})\right)
\end{array}
\label{Definicion VLTau}%
\end{equation}
for any $\pi\in\Pi_{\mathbf{x}}$ and the optimal value function as%
\begin{equation}
V(\mathbf{x})=\sup\nolimits_{\pi\in\Pi_{\mathbf{x}}}V_{\pi}(\mathbf{x}%
).\label{Definicion V}%
\end{equation}
The value $c>0$ is a constant discount factor, and $a_{i}>0$ are the weights
of the dividends paid by the $i$-th company. The integral in
(\ref{Definicion VLTau}) is defined as
\[
\int_{0^{-}}^{t}e^{-cs}\mathbf{a}\cdot d\mathbf{L}_{s}=\mathbf{a}%
\cdot\mathbf{L}_{0}+\int_{0}^{t}e^{-cs}\mathbf{a}\cdot d\mathbf{L}_{s}.
\]
Note that we are allowing to make a lump dividend payment $\mathbf{L}%
_{\overline{\tau}}-\mathbf{L}_{\overline{\tau}^{-}}$ at the switch time
$\overline{\tau}<\tau^{\mathbf{L}}$ and also at time zero.

\begin{remark}
[on the multivariate compound Poisson process]The most important cases of
multivariate compound Poisson process we are considering in the examples
correspond to $m$ independent sources of risk that are coinsured between the
$n$ insurance companies with different proportions. More precisely, let us
assume that there are $m$ independent (univariate) compound Poisson processes
given by%
\begin{equation}
C^{l}(t)=\sum\nolimits_{k=1}^{N_{t}^{l}}u_{k}^{l},\label{CP_Univariadol}%
\end{equation}
where $N_{t}^{l}$ is a Poisson process with intensity $\lambda_{l}$ and
$u_{k}^{l}$ with $k=1,2,...$ are i.i.d. random variables with distribution
$F_{l}$. Assume that the total claim arrival process is given by
\[%
{\displaystyle\sum\nolimits_{j=1}^{N_{t}}}
u_{j}:=\sum\nolimits_{l=1}^{m}C^{l}(t)
\]
and that the $i$-th company pays a proportion $a_{il}$ of any claim of the $l
$-th compound Poisson process $C_{t}^{l}$. We denote $A:=(a_{il})\in
\mathbf{R}^{n\times m}$ with $\sum_{i=1}^{n}a_{il}=1$ and $a_{il}\geq0$. The
compound Poisson process $\sum_{l=1}^{m}C^{l}(t)$ has intensity $\lambda
=\sum_{l=1}^{m}\lambda_{l}$. Furthermore,%
\begin{equation}%
{\displaystyle\sum\nolimits_{k=1}^{N_{t}}}
\mathbf{U}_{k}=A\cdot\left(  C^{1}(t),...,C^{m}(t)\right)  ^{\prime
},\label{CP_Multivariado_RiesgosIndependientes}%
\end{equation}
where $N_{t}=\sum_{l=1}^{m}N_{t}^{l}$ is a compound Poisson process with
intensity $\lambda=\sum_{l=1}^{m}\lambda_{l}$ and multivariate distribution
\[
F(\mathbf{x})=\mathbb{P}(\mathbf{U}\leq\mathbf{x})=\sum\nolimits_{l=1}%
^{m}\frac{\lambda_{l}}{\lambda}F_{l}(\min\nolimits_{1\leq i\leq n,~a_{il}%
\neq0}\{\frac{x_{i}}{a_{il}}\}).
\]

Without loss of generality, we can assume that the columns $\mathbf{a}%
_{l}:=(a_{il})_{i=1,..,n}$ of the matrix $A$ are different, because if
$\mathbf{a}_{l_{1}}=\mathbf{a}_{l_{2}}$, we can regard $C^{l_{1}}(t)+C^{l_{2}%
}(t)$ as just one independent source of risk. For instance, in the special
case in which the $n$ uncontrolled one-dimensional surplus processes of the
companies are independent compound Poisson processes with intensity
$\lambda_{i}$ and claim size distribution $F_{i}(x_{i})$ ($i=1,...,n$), $A$
would be the identity matrix and
\begin{equation}
F(\mathbf{x})=%
{\displaystyle\sum\nolimits_{i=1}^{n}}
\lambda_{i}F_{i}(x_{i})/\lambda\text{.}\label{independent_Surpluses}%
\end{equation}

\end{remark}

\begin{remark}
[on the penalty function $\upsilon$]Consider the multivariate uncontrolled
Poisson process (\ref{UncontrolledSurplusOriginal}) described in the previous
remark. Suppose that the penalty (or reward) function depends on two factors:
(1) which of the $m$ independent compound Poisson processes
(\ref{CP_Univariadol}) make the controlled process exit $\mathbf{R}_{+}^{n}$
and (2) the deficit at this exit time. Let $\mathbf{a}_{l}=(a_{il}%
)_{i=1,..,n}$ be the $l$-th column of $A$, then we have that%
\begin{equation}
\upsilon(\mathbf{x},\mathbf{\alpha})=\sum\nolimits_{l=1}^{m}\upsilon
_{l}(\mathbf{x}-\mathbf{\alpha})I_{\{\mathbf{\alpha}=\beta_{l}\mathbf{a}%
_{l}\text{ with }\beta_{l}>0\}},\label{Funcion NU}%
\end{equation}
where $\upsilon_{l}(\mathbf{X}_{\tau^{\mathbf{L}}}^{\mathbf{L}})$ is the
penalty (or reward) when the process $\mathbf{X}_{t}^{\mathbf{L}}$ exits
$\mathbf{R}_{+}^{n}$ due to a jump of $C^{l}$.

If $n=m=1$, this definition of penalty function $\upsilon$ includes: (1) the
penalty function defined in Gerber and Shiu \cite{GS 1998}, taking
$\upsilon(x,\alpha)=w(x,\left\vert x-\alpha\right\vert )\geq0$; (2) the case
in which the shareholders take care of the deficit at ruin, taking
$\upsilon(x,\alpha)=\alpha-x>0\ $(Dickson and Waters \cite{DicksonWaters2004}%
); (3) the case in which the insurer earns continuously $\Lambda$ as long as
the company is alive. This is equivalent to consider $\upsilon(x,\alpha)=$
$\Lambda/c$ (Thonhauser and Albrecher \cite{TA}).

In the multidimensional framework, the function $\upsilon$ could be negative,
and so considered as a reward. For example, in the case of two companies with
independent compound Poisson processes as in (\ref{independent_Surpluses}), we
can consider the situation in which if one of the companies goes to ruin, the
other survives and continues paying dividends with its own optimal policy. In
this case, $A$ is the $2\times2$ identity matrix and
\begin{equation}
\upsilon(\mathbf{x},\mathbf{\alpha})=-(V_{2}(x_{2})I_{\{x_{1}-\alpha_{1}%
<0\}}+V_{1}(x_{1})I_{\{x_{2}-\alpha_{2}<0\}}%
),\label{Nu the dos companias independientes}%
\end{equation}
where $V_{i}$ is the optimal dividend payment function of the $i$-th company.
Note that $\upsilon(\mathbf{x},\mathbf{\alpha})$ is non-increasing on
$\mathbf{x}$ and non-decreasing on $\mathbf{\alpha}$.
\end{remark}

\begin{remark}
[on the switch-value function $f$]The switch-value function $f(\mathbf{x})$
can be though as the price in which the shareholders can sell the shares when
the controlled current surplus of the $n$ companies is $\mathbf{x}$. It also
can be though as the present value of all the dividends paid in the future
after a change of regime is decided by the manager (this change of regime
could have a cost); for instance, if the manager decides to merge the $n$
companies (that is the $n$ companies put together all their surpluses, pay all
the claims and pay dividends until the joined surplus becomes negative). In
the case of merger,
\begin{equation}
f(\mathbf{x})=V_{M}(x_{1}+x_{2}+...+x_{n}-c_{M})I_{\{x_{1}+x_{2}+...+x_{n}\geq
c_{M}\}}\label{g de merger}%
\end{equation}
where the one-dimensional function $V_{M}$ is the optimal dividend payment
function of the merger of all the companies and $c_{M}\geq0$ is the merger
cost. So, $f$ is right continuous and non decreasing. The case $n=2$, $A$ the
$2\times2$ identity matrix, $\upsilon\ $as in
(\ref{Nu the dos companias independientes}) and
\begin{equation}
f(x_{1},x_{2})=V_{M}(x_{1}+x_{2}-c_{M})I_{\{x_{1}+x_{2}\geq c_{M}%
\}}\label{Merger 2x2}%
\end{equation}
corresponds to the problem of optimal time of merger proposed by Gerber and
Shiu \cite{GS Merger}. The case where no switching is allowed is also included
in this work, just consider $f\ $small enough (see Remark
\ref{V sin obstaculo}).
\end{remark}

In the next proposition we give sufficient conditions under which the function
$V$ is well defined. We say that a function $h:\mathbf{R}_{+}^{n}%
\rightarrow\mathbf{R}$ satisfies the growth condition \textbf{GC} if
\begin{equation}
h(\mathbf{x})/h_{0}(\mathbf{x})\ \text{is upper bounded in }\mathbf{R}_{+}%
^{n}\text{,}\label{gc}%
\end{equation}
where
\begin{equation}
h_{0}(\mathbf{x}):=e^{\frac{c}{2n}\sum_{i=1}^{n}\frac{x_{i}}{p_{i}}%
}.\label{ho}%
\end{equation}

\begin{proposition}
\label{Crecimiento de V} If the functions $f\ $and $S(\mathbf{x}%
):=\sup_{\left\{  \mathbf{\alpha}:\left(  \mathbf{x},\mathbf{\alpha}\right)
\in B\right\}  }\left(  -\upsilon(\mathbf{x},\mathbf{\alpha})\right)  $
satisfy the growth condition \textbf{GC}, then $V$ is well defined, satisfies
the growth condition \textbf{GC} and $V\geq-\mathbb{E}\left(  \left\vert
\upsilon(\mathbf{0},\mathbf{U}_{1})\right\vert \right)  .$
\end{proposition}

\textit{Proof}. Take any initial surplus $\mathbf{x}\geq\mathbf{0}$ and any
admissible strategy $\pi=(\mathbf{L},\overline{\tau})\in\Pi_{\mathbf{x}}$,
since $\mathbf{L}_{t}\leq\mathbf{X}_{t}\leq\mathbf{x}+\mathbf{p}t$, we have
(using integration by parts),%

\[%
\begin{array}
[c]{lll}%
\mathbb{E}_{\mathbf{x}}\left(  \int\nolimits_{0^{-}}^{\tau^{\mathbf{L}}%
\wedge\overline{\tau}}e^{-cs}dL_{i}(s)\right)  & = & \mathbb{E}_{\mathbf{x}%
}\left(  \int\nolimits_{0}^{\tau^{\mathbf{L}}\wedge\overline{\tau}}%
e^{-cs}dL_{i}(s)\right)  +L_{i}(0)\\
& \leq & \mathbb{E}_{\mathbf{x}}\left(  \int\nolimits_{0}^{\tau^{\mathbf{L}%
}\wedge\overline{\tau}}e^{-cs}d(x_{i}+p_{i}s)\right)  +x_{i}\leq x_{i}%
+\frac{p_{i}}{c}.
\end{array}
\]
So%

\[
\mathbb{E}_{\mathbf{x}}\left(  \int_{0^{-}}^{\overline{\tau}\wedge
\tau^{\mathbf{L}}}e^{-cs}\mathbf{a}\cdot d\mathbf{L}_{s}\right)
\leq\mathbf{a}\cdot(\mathbf{x}+\frac{\mathbf{p}}{c})\leq d_{1}e^{\frac{c}%
{2n}\sum_{i=1}^{n}\frac{x_{i}}{p_{i}}}=d_{1}h_{0}(\mathbf{x})
\]
for $d_{1}\geq2n\max\left\{  a_{1},...,a_{n}\right\}  \max\left\{
p_{1},...,p_{n}\right\}  /c$ since%
\[
e^{\frac{c}{2n}\sum_{i=1}^{n}\frac{x_{i}}{p_{i}}}\geq1+\frac{c}{2n}\sum
_{i=1}^{n}\frac{x_{i}}{p_{i}}.
\]
Consider the processes $z(s)=\mathbf{X}_{s}^{\mathbf{L}}{}$ defined in
(\ref{XL}) and let us call $\tau=\tau^{\mathbf{L}}$. We get that
$\mathbf{z}(s)\leq\mathbf{x}+\mathbf{p}s$ and $f$ satisfies (\ref{gc}) in
$\mathbf{R}_{+}^{n}$, so
\[
\mathbb{E}_{\mathbf{x}}\left(  e^{-c\overline{\tau}}f(\mathbf{z}%
_{\overline{\tau}})I_{\{\tau>\overline{\tau}\}}\right)  \leq d_{2}%
e^{\sum_{i=1}^{n}\frac{cx_{i}}{2np_{i}}}=d_{2}h_{0}(\mathbf{x})
\]
for $d_{2}$ large enough. Similarly,
\[
\mathbb{E}_{\mathbf{x}}\left(  -e^{-c\tau}\upsilon(\mathbf{z}_{\tau^{-}%
},\mathbf{z}_{\tau^{-}}-\mathbf{z}_{\tau})I_{\{\tau\leq\overline{\tau}%
\}}\right)  \leq\mathbb{E}_{\mathbf{x}}\left(  e^{-c\tau}I_{\{\tau
\leq\overline{\tau}\}}S(\mathbf{z}_{\tau^{-}})\right)  \leq d_{3}e^{\sum
_{i=1}^{n}\frac{cx_{i}}{2np_{i}}}=d_{3}h_{0}(\mathbf{x})
\]
for $d_{3}$ large enough. Then $V_{\pi}$ (and so $V$) satisfy the growth
condition (\ref{gc}) in $\mathbf{R}_{+}^{n}$. Finally, since $\tau$ is the
first time that the controlled process $\mathbf{X}^{\mathbf{L}}$ leaves
$\mathbf{R}_{+}^{n}$, calling $\mathbf{U}_{k_{0}}$ the jump size at $\tau$, we
have $\mathbf{z}_{\tau}=\mathbf{z}_{\tau^{-}}-\mathbf{U}_{k_{0}}\geq
\mathbf{0}-\mathbf{U}_{k_{0}}$. Since $-\upsilon(\mathbf{x},\mathbf{\alpha})$
is non-decreasing on $\mathbf{x}$ and non-increasing on $\mathbf{\alpha}$, we
obtain taking the strategy with no switching and no dividend payment, that
\[
V(\mathbf{x})\geq\mathbb{E}_{\mathbf{x}}\left(  -e^{-c\tau}\upsilon
(\mathbf{z}_{\tau^{-}},\mathbf{U}_{k_{0}}\right)  \geq\mathbb{E}_{\mathbf{x}%
}\left(  -e^{-c\tau}\upsilon(\mathbf{0},\mathbf{U}_{k_{0}})\right)
\geq-\mathbb{E}\left(  \left\vert \upsilon(\mathbf{0},\mathbf{U}%
_{1})\right\vert \right)  \text{. }\blacksquare
\]

\begin{remark}
Let us extend the definition of $\upsilon$ to the closure of $B$ as
$\upsilon(\mathbf{x},\mathbf{\alpha})=\inf_{\mathbf{\beta}\geq\mathbf{\alpha
},\left(  \mathbf{x},\mathbf{\beta}\right)  \in B}\upsilon(\mathbf{x}%
,\mathbf{\beta})$. Since $-\upsilon(\mathbf{x},\mathbf{\alpha})$ is
non-decreasing on $\mathbf{x}$ and non-increasing on $\mathbf{\alpha}$ then
\[
\sup\nolimits_{\mathbf{\alpha}\geq\mathbf{0},\mathbf{x}-\mathbf{\alpha}%
\notin\mathbf{R}_{+}^{n}}(-\upsilon(\mathbf{x},\mathbf{\alpha}))\leq\max
{}_{i=1,...n}(-\upsilon(\mathbf{x},x_{i}\mathbf{e}_{i}))
\]
and so the assumption on $\upsilon$ of Proposition \ref{Crecimiento de V}
becomes that $\max{}_{i=1,...n}\left(  -\upsilon(\mathbf{x},x_{i}%
\mathbf{e}_{i})\right)  $ satisfies the growth condition \textbf{GC}.
\end{remark}

\begin{remark}
\label{V sin obstaculo} By Proposition \ref{Crecimiento de V}, taking any
switch-value function $f<-\mathbb{E}\left(  \left\vert \upsilon(\mathbf{0}%
,\mathbf{U}_{1})\right\vert \right)  $, it is never optimal to switch. So, the
problem of maximizing the expected cumulative discounted dividend payments
until $\tau^{\mathbf{L}}$ (without the possibility of switching) is a
particular case of the problem (\ref{Definicion V}).
\end{remark}

\begin{remark}
\label{Monotonia de Obstaculos} Consider $f_{1}\leq f_{2}$ and $\upsilon
_{1}\geq\upsilon_{2}$ . Let $V_{f_{1},\upsilon_{1}}$ and $V_{f_{2}%
,\upsilon_{2}}$ be the corresponding optimal value functions, then it is
straightforward to see that $V_{f_{1},\upsilon_{1}}\leq V_{f_{2},\upsilon_{2}%
}$.
\end{remark}

\begin{remark}
Since the optimal dividend payment function in the one-dimensional problem has
linear growth, see for instance Proposition 1.2 in \cite{AM Libro}; the
functions (\ref{Nu the dos companias independientes})\ and (\ref{g de merger})
satisfy the conditions of Proposition \ref{Crecimiento de V}\textbf{.}
\end{remark}

In the next proposition, we show that $V$ is increasing and locally Lipschitz
(so it is absolutely continuous).

\begin{proposition}
\label{V Lipschitz} $V$ is increasing, locally Lipschitz in $\mathbf{R}%
_{+}^{n}$ and satisfies for each $\mathbf{x}\in$ $\mathbf{R}_{+}^{n}$, $h>0$
and $1\leq i\leq n$,
\[
a_{i}h\leq V(\mathbf{x}+h\mathbf{e}_{i})-V(\mathbf{x})\leq(e^{(c+\lambda
)h/p_{i}}-1)V(\mathbf{x}).
\]

\end{proposition}

\textit{Proof}. Given $h>0$ and $\mathbf{x}\in$ $\mathbf{R}_{+}^{n}$, consider
for each $\varepsilon>0$ an admissible strategy $\pi_{\mathbf{x}}=\left(
\mathbf{L},\overline{\tau}\right)  \in\Pi_{\mathbf{x}}$ such that $V_{\pi
}(\mathbf{x})\geq V(\mathbf{x})-\varepsilon$. \ Let us define an strategy
$\widetilde{\pi}\in\Pi_{\mathbf{x}+h\mathbf{e}_{i}}$ as follows: the $i$-th
company pays immediately $h$ as dividends and then follows the strategy $\pi$
$\in\Pi_{\mathbf{x}}$. For each $\varepsilon>0$, we get
\[
V(\mathbf{x}+h\mathbf{e}_{i})\geq V_{\widetilde{\pi}}(\mathbf{x}%
+h\mathbf{e}_{i})=V_{\pi}(\mathbf{x})+a_{i}h\geq V(\mathbf{x})-\varepsilon
+a_{i}h,
\]
so we obtain the first inequality.

Now consider for each $\varepsilon>0$ and $1\leq i\leq n$, a strategy
$\pi=(\mathbf{L},\overline{\tau})\in\Pi_{\mathbf{x}+h\mathbf{e}_{i}}$ such
that
\[
V_{\pi}(\mathbf{x}+h\mathbf{e}_{i})\geq V(\mathbf{x}+h\mathbf{e}%
_{i})-\varepsilon.
\]
Take now the following admissible strategy $\widetilde{\pi}=(\widetilde
{\mathbf{L}},\widetilde{\tau})\in\Pi_{\mathbf{x}}$ starting with surplus
$\mathbf{x}$: the $i$-th company pays no dividends and the other companies pay
all the incoming premium as dividends as long as $\mathbf{X}_{t}%
^{\widetilde{\mathbf{L}}}<\mathbf{x}+h\mathbf{e}_{i}$; after the current
surplus reaches $\mathbf{x}+h\mathbf{e}_{i}$, follow strategy $\pi$. Let us
call $\tau^{\widetilde{L}}$ the exit time of the process $\mathbf{X}%
_{t}^{\widetilde{\mathbf{L}}}$. If
\[
\tau:=\min\{t:\mathbf{X}_{t}^{\widetilde{\mathbf{L}}}=\mathbf{x}%
+h\mathbf{e}_{i}\},
\]
then $\widetilde{\tau}=\tau+\overline{\tau}$ and we get that $p_{i}\tau\geq
h$. So,%

\[%
\begin{array}
[c]{lll}%
V(\mathbf{x}) & \geq & V_{\widetilde{\pi}}(\mathbf{x})\geq V_{\pi}%
(\mathbf{x}+h\mathbf{e}_{i})\mathbb{E}\left(  e^{-c\frac{h}{p_{i}}}%
I_{\{\tau<\tau^{\widetilde{L}}\}}\right)  \geq\left(  V(\mathbf{x}%
+h\mathbf{e}_{i})-\varepsilon\right)  e^{-c\frac{h}{p_{i}}}\mathbb{P}%
(\tau<\tau^{\widetilde{L}})\\
& \geq & \left(  V(\mathbf{x}+h\mathbf{e}_{i})-\varepsilon\right)
e^{-c\frac{h}{p_{i}}}\mathbb{P}(\tau_{1}>\frac{h}{p_{i}})=\left(
V(\mathbf{x}+h\mathbf{e}_{i})-\varepsilon\right)  e^{-\left(  c+\lambda
\right)  \frac{h}{p_{i}}}\text{,}%
\end{array}
\]
where $\tau_{1}$ is the time of the first jump; so we get the second
inequality. $\blacksquare$

In order to distinguish the jumps of the controlled process due to the jumps
of the uncontrolled process from the ones due to lump dividend payments, let
us define an auxiliary process which includes the jump of the uncontrolled
process occurring at time $t$ but excludes the lump dividend payment occurring
at this time as
\begin{equation}
\mathbf{\check{X}}_{t}^{\mathbf{L}}=\mathbf{X}_{t}-\mathbf{L}_{t^{-}%
}=\mathbf{X}_{t^{-}}^{\mathbf{L}}-\left(  \mathbf{X}_{t^{-}}-\mathbf{X}%
_{t}\right)  \text{.}\label{Xv}%
\end{equation}
Note that $\mathbf{\check{X}}_{t}^{\mathbf{L}}=\mathbf{X}_{t^{-}}^{\mathbf{L}%
}-\mathbf{U}_{k}$ if $t=\tau_{k}$ and $\mathbf{\check{X}}_{t}^{\mathbf{L}%
}=\mathbf{X}_{t^{-}}^{\mathbf{L}}$ otherwise. Also, $\mathbf{X}_{\tau
^{\mathbf{L}}}^{\mathbf{L}}=\mathbf{\check{X}}_{\tau^{\mathbf{L}}}%
^{\mathbf{L}}$ because no dividends are paid at the exit time $\tau
^{\mathbf{L}}$.

\section{HJB equation\label{Seccion Viscosidad}}

In this section we show that the optimal value function $V$ defined in
(\ref{Definicion V}) is a viscosity solution of the corresponding
Hamilton-Jacobi-Bellman (HJB) equation; moreover we characterize the optimal
value function as the smallest viscosity supersolution with growth condition
\textbf{GC}. We also give a verification result for $V$. These results are a
generalization to the multidimensional case of the ones given in Section 3 of
\cite{AM Switching} for the one dimensional case.

The HJB equation of problem (\ref{Definicion V}) can be written as%
\begin{equation}
\max\{\mathbf{a}-\nabla V(\mathbf{x}),\mathcal{L}(V)(\mathbf{x}),f(\mathbf{x}%
)-V(\mathbf{x})\}=0,\label{HJB}%
\end{equation}
where%
\begin{equation}
\mathcal{L}(V)(\mathbf{x})=\mathbf{p\cdot}\nabla V(\mathbf{x})-(c+\lambda
)V(\mathbf{x})+\mathcal{I}(V)(\mathbf{x})-\mathcal{R}(\mathbf{x}%
),\label{L para Cramer-Lundberg}%
\end{equation}

\begin{equation}
\mathcal{I}(W)(\mathbf{x}):=\lambda\int_{\left(  \mathbf{x}-\mathbf{\alpha
}\right)  \in\mathbf{R}_{+}^{n}}W(\mathbf{x}-\mathbf{\alpha})dF(\mathbf{\alpha
})\text{ and }\mathcal{R}(\mathbf{x}):=\lambda\int_{\left(  \mathbf{x}%
-\mathbf{\alpha}\right)  \notin\mathbf{R}_{+}^{n}}\upsilon(\mathbf{x}%
,\mathbf{\alpha})dF(\mathbf{\alpha})\text{.}\label{Def R e I}%
\end{equation}

As usual, the operator $\mathcal{L}$ is the discounted infinitesimal generator
of the uncontrolled surplus process $\mathbf{X}_{t}$ defined in
(\ref{UncontrolledSurplusOriginal}); that is, for any continuously
differentiable function $W:\mathbf{R}_{+}^{n}\rightarrow\mathbf{R}$\textbf{,}
we have
\[
\mathcal{L}(W)(\mathbf{x})=\lim_{t\searrow0}\frac{\mathbb{E}_{\mathbf{x}%
}\left(  e^{-ct}W(\mathbf{X}_{t})-W(\mathbf{x})\right)  }{t}.
\]
Thus, if $W$ is a solution of $\mathcal{L}(W)=0$ in an open set in
$\mathbf{R}_{+}^{n}$, then the process $e^{-ct}W(\mathbf{X}_{t})$ is a
martingale in this set.

The HJB equation implies that $\mathcal{L}(V)\leq0$, the condition
$\mathcal{L}(V)=0$ in an open set in $\mathbf{R}_{+}^{n}$ would suggest that
(locally) the optimal dividend strategy consists on paying no dividends as
long as the current surplus is in this set. The HJB equation also implies that
$V$ is always above $f$, so $f$ can be interpreted as \textit{an obstacle }in
equation (\ref{HJB}). Moreover, the condition $V_{x_{i}}(\mathbf{x})=a_{i}$ in
an open set means that (locally) the optimal dividend strategy should be the
one in which the $i$-th company pays immediately a lump sum as dividends.

We prove in this section that, under the assumption%
\begin{equation}
\mathcal{R}:\mathbf{R}_{+}^{n}\rightarrow\mathbf{R}\text{ is continuous,}%
\label{R connttinua}%
\end{equation}
the value function $V$ is a viscosity solution of the HJB equation
(\ref{HJB}). From now on, we assume that this assumption holds.

Crandall and Lions \cite{CL} introduced the concept of viscosity solutions for
first-order Hamilton-Jacobi equations. It is the standard tool for studying
HJB equations, see for instance Fleming and Soner \cite{FS}. In the context of
using viscosity solutions for the problem of dividend payment optimization in
the one-dimensional case, see for instance \cite{AM Libro}.

\begin{definition}
\label{NuevaDefinicionSubySuper}A locally Lipschitz function $\underline
{u}:\mathbf{R}_{+}^{n}\rightarrow\mathbf{R}$\ is a viscosity subsolution of
(\ref{HJB}) at $\mathbf{x}\in$ $\mathbf{R}_{+}^{n}$\ if any continuously
differentiable function $\psi\ $defined in $\mathbf{R}_{+}^{n}$\ with
$\psi(\mathbf{x})=$\ $\underline{u}(\mathbf{x})$\ such that $\underline
{u}-\psi$\ reaches the maximum at $\mathbf{x}$\ satisfies
\[
\max\{\mathbf{a}-\nabla\psi(\mathbf{x}),\mathcal{L}(\psi)(\mathbf{x}%
),f(\mathbf{x})-\psi(\mathbf{x})\}\geq0,
\]
and a locally Lipschitz function $\overline{u}$\ $:\mathbf{R}_{+}%
^{n}\rightarrow\mathbf{R}$\ is a viscosity supersolution of (\ref{HJB}) at
$\mathbf{x}\in$ $\mathbf{R}_{+}^{n}$\ if any continuously differentiable
function $\varphi$ defined in $\mathbf{R}_{+}^{n}$\ \ with $\varphi
(\mathbf{x})=$\ $\overline{u}(\mathbf{x})$\ such that $\overline{u}-\varphi
$\ reaches the minimum at $\mathbf{x}$\ satisfies
\[
\max\{\mathbf{a}-\nabla\varphi(\mathbf{x}),\mathcal{L}(\varphi)(\mathbf{x}%
),f(\mathbf{x})-\varphi(\mathbf{x})\}\leq0.
\]
Finally, a locally Lipschitz function $u$\ $:\mathbf{R}_{+}^{n}\rightarrow
\mathbf{R}$\ is a viscosity solution of (\ref{HJB}) if it is both a viscosity
subsolution and a viscosity supersolution at any $\mathbf{x}\in\mathbf{R}%
_{+}^{n}$.
\end{definition}

In order to prove that $V$ is a viscosity solution of the HJB equation we need
to use the following two lemmas. The first one states the Dynamic Programming
Principle (DPP), its proof follows from standard arguments, see for instance
Lemma 1.2 of \cite{AM Libro}. The proof of the second one is in the Appendix.

\begin{lemma}
\label{DPP} Given any $\mathbf{x}\in\mathbf{R}_{+}^{n}$\ and any finite
stopping time $\widetilde{\tau}$,\ we have that the function $V$ defined in
(\ref{Definicion V}) satisfies $V(\mathbf{x})=\sup_{\pi=(L,\overline{\tau}%
)\in\Pi_{\mathbf{x}}}v_{\pi,\widetilde{\tau}}(\mathbf{x})$, where
\[%
\begin{array}
[c]{lll}%
v_{\pi,\widetilde{\tau}}(\mathbf{x}) & = & \mathbb{E}_{\mathbf{x}}\left(
\int\nolimits_{0^{-}}^{\overline{\tau}\wedge\tau^{\mathbf{L}}\wedge
\widetilde{\tau}}e^{-cs}\mathbf{a}\cdot d\mathbf{L}_{s}+e^{-c(\overline{\tau
}\wedge\tau^{\mathbf{L}}\wedge\widetilde{\tau})}(I_{\{\overline{\tau}%
\wedge\widetilde{\tau}<\tau^{\mathbf{L}}\}}V(\mathbf{X}_{\overline{\tau}%
\wedge\widetilde{\tau}}^{\mathbf{L}})\right) \\
&  & -\mathbb{E}_{\mathbf{x}}\left(  I_{\{\tau^{\mathbf{L}}\leq\overline{\tau
}\wedge\widetilde{\tau}\}}\upsilon(\mathbf{X}_{\tau^{\mathbf{L}}{}^{-}%
}^{\mathbf{L}}{},\mathbf{X}_{\tau^{\mathbf{L}}{}^{-}}^{\mathbf{L}}%
-\mathbf{X}_{\tau^{\mathbf{L}}}^{\mathbf{L}})\right)  .
\end{array}
\]

\end{lemma}

\begin{lemma}
\label{Dynkins} Given any continuously differentiable function $g:\mathbf{R}%
_{+}^{n}\rightarrow\mathbf{R}$, any admissible strategy $\pi=(\mathbf{L}%
,\overline{\tau})\in\Pi_{\mathbf{x}}$ and any finite stopping time $\tau
\leq\tau^{\mathbf{L}}$, consider
\[
\mathbf{L}_{t}=\int\nolimits_{0}^{t}d\mathbf{L}_{s}^{c}+\sum\nolimits_{0\leq
s\leq t}\Delta\mathbf{L}_{s},
\]
where $\Delta\mathbf{L}_{s}=\mathbf{L}_{s}-\mathbf{L}_{s^{-}}$ and
$\mathbf{L}_{s}^{c}$ is a continuous and non-decreasing process. Then we have%
\[%
\begin{array}
[c]{l}%
(g(\mathbf{X}_{\tau}^{\mathbf{L}})I_{\{\tau<\tau^{\mathbf{L}}\}}%
-\upsilon(\mathbf{X}_{\tau^{-}}^{\mathbf{L}},\mathbf{X}_{\tau^{-}}%
^{\mathbf{L}}-\mathbf{X}_{\tau}^{\mathbf{L}})I_{\{\tau=\tau^{\mathbf{L}}%
\}})e^{-c\tau}-g(\mathbf{x})\\%
\begin{array}
[c]{ll}%
= & \int\nolimits_{0}^{\tau}\mathcal{L}(g)(\mathbf{X}_{s^{-}}^{\mathbf{L}%
})e^{-cs}ds-\int_{0^{-}}^{\tau}e^{-cs}\mathbf{a}\cdot d\mathbf{L}_{s}%
+\int\nolimits_{0}^{\tau}e^{-cs}(\mathbf{a}-\nabla g(\mathbf{X}_{s^{-}%
}^{\mathbf{L}}))\mathbf{\cdot}d\mathbf{L}_{s}^{c}\\
& +\sum\limits_{\mathbf{L}_{s}\neq\mathbf{L}_{s^{-}},s\leq\tau}e^{-cs}%
\int\nolimits_{0}^{1}(\mathbf{a}-\nabla g(\mathbf{\check{X}}_{s}^{\mathbf{L}%
}-\gamma\Delta\mathbf{L}_{s})\mathbf{\cdot}\Delta\mathbf{L}_{s})d\gamma
+M(\tau);
\end{array}
\end{array}
\]
where $M(t)$ is a martingale with zero expectation.
\end{lemma}

\begin{proposition}
\label{Prop V is a viscosity solution}The optimal value function $V$ is a
viscosity solution of the HJB equation (\ref{HJB}) at any $\mathbf{x}$ in the
interior of $\mathbf{R}_{+}^{n}$.
\end{proposition}

\textit{Proof}. Let us show that $V$ is a viscosity supersolution at any
$\mathbf{x}$ in the interior of $\mathbf{R}_{+}^{n}$. The inequality $V\geq f$
follows from the definition (\ref{Definicion V}) taking $\overline{\tau}=0$.
Given any initial surplus $\mathbf{x}$ in the interior of $\mathbf{R}_{+}^{n}$
and any $\mathbf{l}\in\mathbf{R}_{+}^{n}$, take $h>0$ small enough such that
$h(\mathbf{l}-\mathbf{p})<\mathbf{x}$\textbf{.} Consider the dividend payment
strategy $\mathbf{L}_{t}=\mathbf{l}t$ for $t<h\wedge$ $\tau_{1}$ and
$\mathbf{L}_{t}=\mathbf{l}(h\wedge\tau_{1})$ for $t\geq h\wedge$ $\tau_{1}$;
also consider a switch time $\overline{\tau}>\tau^{\mathbf{L}}$. Using Lemma
\ref{DPP} with stopping time $\widetilde{\tau}=h\wedge\tau_{1}$, we get%
\[
V(\mathbf{x})\geq\mathbb{E}_{\mathbf{x}}\left(  \mathbf{a}\cdot\mathbf{l}%
\int\nolimits_{0}^{h\wedge\tau_{1}}e^{-cs}ds+e^{-c\left(  h\wedge\tau
_{1}\right)  }(I_{\{h\wedge\tau_{1}<\tau^{\mathbf{L}}\}}V(\mathbf{X}%
_{h\wedge\tau_{1}}^{\mathbf{L}})-I_{\{\tau^{\mathbf{L}}=\tau_{1}\leq
h\}}\upsilon(\mathbf{X}_{\tau_{1}^{-}}^{\mathbf{L}},\mathbf{U}_{1}))\right)  .
\]
Let $\varphi$ be a test function for supersolution of (\ref{HJB}) at
$\mathbf{x}$ as in Definition \ref{NuevaDefinicionSubySuper}. We have,%
\[%
\begin{array}
[c]{lll}%
\varphi(\mathbf{x}) & = & V(\mathbf{x})\\
& \geq & \mathbb{E}_{\mathbf{x}}\left(  \mathbf{a}\cdot\mathbf{l}%
\int\nolimits_{0}^{h\wedge\tau_{1}}e^{-cs}ds\right) \\
&  & +\mathbb{E}_{\mathbf{x}}\left(  e^{-c\left(  h\wedge\tau_{1}\right)
}(I_{\{h\wedge\tau_{1}<\tau^{\mathbf{L}}\}}\varphi(\mathbf{X}_{h\wedge\tau
_{1}}^{\mathbf{L}})-I_{\{\tau^{\mathbf{L}}=\tau_{1}\leq h\}}\upsilon
(\mathbf{X}_{\tau_{1}^{-}}^{\mathbf{L}},\mathbf{U}_{1}))\right)  .
\end{array}
\]
We can write%
\[%
\begin{array}
[c]{l}%
\mathbb{E}_{\mathbf{x}}\left(  e^{-c\left(  h\wedge\tau_{1}\right)
}(I_{\{h\wedge\tau_{1}<\tau^{\mathbf{L}}\}}\varphi(\mathbf{X}_{h\wedge\tau
_{1}}^{\mathbf{L}})-I_{\{\tau^{\mathbf{L}}=\tau_{1}\leq h\}}\upsilon
(\mathbf{X}_{\tau_{1}^{-}}^{\mathbf{L}},\mathbf{U}_{1}))\right) \\%
\begin{array}
[c]{ll}%
= & \mathbb{E}_{\mathbf{x}}\left(  I_{\{h<\tau_{1}\text{ }\}}e^{-ch}%
\varphi(\mathbf{x}+\left(  \mathbf{p}-\mathbf{l}\right)  h\right) \\
& +\mathbb{E}_{\mathbf{x}}\left(  I_{\{\tau_{1}\leq h\}}I_{\{\mathbf{U}%
_{1}\leq\mathbf{x}+\left(  \mathbf{p}-\mathbf{l}\right)  \tau_{1}\}}%
e^{-c\tau_{1}}\varphi(\mathbf{x}+\left(  \mathbf{p}-\mathbf{l}\right)
\tau_{1}-\mathbf{U}_{1})\right) \\
& -\mathbb{E}_{\mathbf{x}}\left(  I_{\{\tau_{1}\leq h\}}I_{\{\mathbf{U}%
_{1}\nleqslant\mathbf{x}+\left(  \mathbf{p}-\mathbf{l}\right)  \tau_{1}%
\}}e^{-c\tau_{1}}\upsilon(\mathbf{x}+\left(  \mathbf{p}-\mathbf{l}\right)
\tau_{1},\mathbf{U}_{1})\right)  .
\end{array}
\end{array}
\]
Therefore, using that $\mathcal{R}$ is continuous,%
\begin{equation}%
\begin{array}
[c]{lll}%
0 & \geq & (\mathbf{a}\cdot\mathbf{l})\lim_{h\rightarrow0^{+}}\frac{1}%
{h}\mathbb{E}_{\mathbf{x}}\left(  \int\nolimits_{0}^{h\wedge\tau_{1}}%
e^{-cs}ds\right)  +\lim_{h\rightarrow0^{+}}\frac{1}{h}\left(  e^{-(\lambda
+c)h}\varphi(\mathbf{x}+\left(  \mathbf{p}-\mathbf{l}\right)  h)-\varphi
(\mathbf{x})\right) \\
&  & +\lim_{h\rightarrow0^{+}}\frac{1}{h}\mathbb{E}_{\mathbf{x}}\left(
I_{\{\tau_{1}\leq h\}}I_{\{\mathbf{U}_{1}\leq\mathbf{x}+\left(  \mathbf{p}%
-\mathbf{l}\right)  \tau_{1}\}}e^{-c\tau_{1}}\varphi(\mathbf{x}+\left(
\mathbf{p}-\mathbf{l}\right)  \tau_{1}-\mathbf{U}_{1})\right) \\
&  & -\lim_{h\rightarrow0^{+}}\frac{1}{h}\mathbb{E}_{\mathbf{x}}\left(
I_{\{\tau_{1}\leq h\}}I_{\{\mathbf{U}_{1}\leq\mathbf{x}+\left(  \mathbf{p}%
-\mathbf{l}\right)  \tau_{1}\}}e^{-c\tau_{1}}\upsilon(\mathbf{x}+\left(
\mathbf{p}-\mathbf{l}\right)  \tau_{1},\mathbf{U}_{1})\right) \\
& = & \mathbf{a}\cdot\mathbf{l}-\left(  c+\lambda\right)  \varphi
(\mathbf{x})+\left(  \mathbf{p}-\mathbf{l}\right)  \mathbf{\cdot}\nabla
\varphi(\mathbf{x})+\mathcal{I}(\varphi)(\mathbf{x})-\mathcal{R}(\mathbf{x}).
\end{array}
\nonumber
\end{equation}
And so $\mathcal{L}(\varphi)(\mathbf{x})+\mathbf{l}$\textbf{$\cdot$}$\left(
\mathbf{a}-\nabla\varphi(\mathbf{x})\right)  \leq0$. Taking $\mathbf{l}%
=\mathbf{0}$, we get $\mathcal{L}(\varphi)(\mathbf{x})\leq0$; taking
$\mathbf{l}=l\mathbf{e}_{i}$ with $l\rightarrow\infty$ ($1\leq i\leq n$), we
obtain $\mathbf{a}-\nabla\varphi(\mathbf{x})\leq\mathbf{0}$. So $V$ is a
viscosity supersolution at the point $\mathbf{x}$.

We omit the proof that $V$ is a viscosity subsolution in the interior of
$\mathbf{R}_{+}^{n}$. This result follows from Lemma \ref{Dynkins} and the
proof is similar to the ones of Proposition 3.2 in \cite{AM Switching} for the
unidimensional case with switching and of Proposition 3.2 in \cite{AlAZMU} for
the multidimensional case with no switching. $\blacksquare$

\begin{remark}
\label{Muchas soluciones viscosas} In general, we cannot expect to have
uniqueness of viscosity solutions of the HJB equation (\ref{HJB}). Take for
instance the two dimensional case with independent companies, the switch
function $f$ given in (\ref{Merger 2x2}) and the function $\upsilon$ defined
in (\ref{Nu the dos companias independientes}). Consider the function
$W_{k}(\mathbf{x}):=x_{1}+x_{2}+k$ for $\mathbf{x\in R}_{+}^{2},$ and take
$k_{0}$ large enough such that,%
\[
k_{0}>\frac{p_{1}+p_{2}}{c}\text{, }V_{2}(z)<z+k_{0},V_{1}(z)<z+k_{0}\text{
and }V_{M}(z)<z+c_{M}+k_{0}%
\]

for $z\geq0$. Hence, we have for all $k\geq k_{0}$ that $\nabla W_{k}%
(\mathbf{x})-\mathbf{1=0}$; $f(\mathbf{x})-W_{k}(\mathbf{x})<0$ and%
\[%
\begin{array}
[c]{lll}%
\mathcal{L}(W_{k})(\mathbf{x}) & \leq & p_{1}+p_{2}-\left(  c+\lambda\right)
W_{k}(\mathbf{x})+\lambda_{1}W_{k}(\mathbf{x})F_{1}(x_{1})+\lambda_{2}%
W_{k}(\mathbf{x})F_{2}(x_{2})\\
&  & +\lambda_{1}(x_{2}+k)(1-F_{1}(x_{1}))+\lambda_{2}(x_{1}+k\text{
})(1-F_{2}(x_{2}))\\
& \leq & p_{1}+p_{2}-\left(  c+\lambda\right)  W_{k}(\mathbf{x})+\lambda
_{1}W_{k}(\mathbf{x})+\lambda_{2}W_{k}(\mathbf{x})\\
& = & -c\left(  x_{1}+x_{2}+K\right)  +p_{1}+p_{2}\\
& < & 0.
\end{array}
\]

Therefore, there are infinitely many viscosity solutions of the HJB equation
(\ref{HJB}).
\end{remark}

The following lemma states that any viscosity supersolution with the
appropriate growth condition is above the value function of a family of
admissible strategies. The proof is in the Appendix.

\begin{lemma}
\smallskip\label{SupersolucionMayor-ValueFunction} Fix $\mathbf{x}%
\in\mathbf{R}_{+}^{n}$, let $\overline{u}$\ be a viscosity supersolution of
(\ref{HJB}) satisfying growth condition (\ref{gc}) and take any admissible
strategy $\pi=(\mathbf{L},\overline{\tau})\in\Pi_{\mathbf{x}}$, then
$\overline{u}(\mathbf{x})\geq V_{\pi}(\mathbf{x})$ (and so $\overline
{u}(\mathbf{x})\geq V(\mathbf{x})$).
\end{lemma}

There are not natural boundary conditions for the optimal value function $V$
(see for instance Section 1.6 of \cite{AM Libro} for a discussion about it in
the one dimensional case). As a consequence of the previous lemma, we get the
following characterization result:

\begin{theorem}
\label{characterization} The optimal value function $V$ can be characterized
as the smallest viscosity supersolution of the HJB equation (\ref{HJB})
satisfying growth condition (\ref{gc}).
\end{theorem}

Also, we obtain immediately the next verification theorem.

\begin{theorem}
\label{verification result} Consider a family of strategies $\left(
\pi_{\mathbf{x}}\right)  _{\mathbf{x}\in\mathbf{R}_{+}^{n}}$ where each
$\pi_{\mathbf{x}}\in\Pi_{\mathbf{x}}$. If the function $W(\mathbf{x}%
)=V_{\pi_{\mathbf{x}}}(\mathbf{x})$ is a viscosity supersolution of the HJB
equation (\ref{HJB}) in the interior of $\mathbf{R}_{+}^{n}$, then $W$ is the
optimal value function. Also, if for each $k\geq1$ there exists a family of
strategies $\left(  \pi_{\mathbf{x}}^{k}\right)  _{\mathbf{x}\in\mathbf{R}%
_{+}^{n}}$ with $\pi_{\mathbf{x}}^{k}\in\Pi_{\mathbf{x}}$ such that
$W(\mathbf{x}):=\lim_{k\rightarrow\infty}V_{\pi_{\mathbf{x}}^{k}}(\mathbf{x})
$ is a viscosity supersolution of the HJB equation (\ref{HJB}) in the interior
of $\mathbf{R}_{+}^{n}$, then $W$ is the optimal value function.
\end{theorem}

\begin{remark}
It is easy to show that the function $h_{0}$ defined\ in (\ref{ho}) satisfies
that $bh_{0}$ are supersolutions of the HJB equation (\ref{HJB}) for all $b$
large enough.
\end{remark}

\section{Discrete Scheme\label{Discrete Scheme}}

In this section we construct a family of admissible strategies for any point
in a suitable grid and then extend it to $\mathbf{R}_{+}^{n}$. We will show in
the next section that the value function of these strategies converge to the
optimal value function $V$ as the mesh size goes to zero.

Given any approximation parameter $\delta>0$, we define the grid domain
\[
\mathcal{G}^{\delta}:=\left\{  (m_{1}\delta p_{1},...,m_{n}\delta
p_{n}):\mathbf{m}\in\mathbf{N}_{0}^{n}\right\}  \text{.}%
\]
The idea of the numerical scheme is to find, at each point of the grid
$\mathcal{G}^{\delta}$, the best local strategy among the ones suggested by
the operators of the HJB equation (\ref{HJB}); these possible local strategies
are: none of the companies pay dividends, one of the companies pays a lump sum
as dividends, or the manager of the company opts to switch immediately. We
modify these local strategies in such a way that the controlled surplus lies
in the grid after the arrival of a jump of the uncontrolled process. In order
to do that, let us introduce the functions $g^{\delta}:\mathbf{N}_{0}%
^{n}\rightarrow\mathbf{R}_{+}^{n}$ \ which relates the indices with the
corresponding points of the grid and $\rho^{\delta}:\mathbf{R}_{+}%
^{n}\rightarrow\mathbf{N}_{0}^{n}$ which assigns to each point points
$\mathbf{x}$ in $\mathbf{R}_{+}^{n}$ the index of the closest point of the
grid below $\mathbf{x}$. More precisely,
\[
g^{\delta}(\mathbf{m})=(p_{1}\delta m_{1},...,p_{n}\delta m_{n})~\text{and
}\rho^{\delta}(\mathbf{x}):=\max\{\mathbf{m}\in\mathbf{N}_{0}^{n}:g^{\delta
}(\mathbf{m})\leq\mathbf{x}\}\text{;}%
\]
we can also write
\[
\rho^{\delta}(\mathbf{x})=(\left[  \frac{x_{1}}{\delta p_{1}}\right]
,...,\left[  \frac{x_{n}}{\delta p_{n}}\right]  )\in\mathbf{N}_{0}^{n}%
\]
where $\left[  .\right]  $ means the integer part in each coordinate. Note
that $\rho^{\delta}\ $is the left-inverse function of $g^{\delta}$ and that%
\[
\left\langle \mathbf{x}\right\rangle ^{\delta}:=g^{\delta}(\rho^{\delta
}(\mathbf{x}))=\max\{\mathbf{y}\in\mathcal{G}^{\delta}:\mathbf{y}%
\leq\mathbf{x}\}\text{.}%
\]

Given any current surplus $g^{\delta}(\mathbf{m})\in\mathcal{G}^{\delta}$, let
$\tau$ and $\mathbf{U}$ be the arrival time and the size of the next jump of
the uncontrolled process. We first define the $n+2$ possible control actions
at any point of the grid $\mathcal{G}^{\delta}$ as follows.

\begin{itemize}
\item Control action $\mathbf{E}_{0}$: Pay no dividends up to the time
$\delta\wedge\tau$. In the case that $\delta<\tau$, the uncontrolled surplus
at time $\delta$ is $g^{\delta}\left(  \mathbf{m}+\mathbf{1}\right)
\in\mathcal{G}^{\delta}$; and if $\delta\geq\tau$, the uncontrolled surplus at
time $\tau$ is
\[
g^{\delta}(\mathbf{m})+\tau\mathbf{p}-\mathbf{U}.
\]
If this vector is in $\mathbf{R}_{+}^{n}$, the companies pay immediately the
minimum amount of dividends in such a way that the controlled surplus lies in
a point of the grid; this end surplus can be written as $g^{\delta}\left(
\mathbf{k}\right)  $, where
\[
\mathbf{k}=\rho^{\delta}(g^{\delta}(\mathbf{m})+\tau\mathbf{p}-\mathbf{U}).
\]
The amount paid as dividends is equal to
\[
g^{\delta}\left(  \mathbf{m}-\mathbf{k}\right)  +\tau\mathbf{p}-\mathbf{U}.
\]
In the case that the surplus $g^{\delta}(\mathbf{m})+\tau\mathbf{p}%
-\mathbf{U}\notin\mathbf{R}_{+}^{n}$ at time $\tau\leq$ $\delta$, the process stops.

\item Control actions $\mathbf{E}_{i}$ with $i=1,...,n$: The $i$-th company
pays immediately $p_{i}\delta$ as dividends, so the controlled surplus becomes
$g^{\delta}\left(  \mathbf{m}-\mathbf{e}_{i}\right)  \in\mathcal{G}^{\delta}$.
The control action $\mathbf{E}_{i}$\textbf{\ }can only be applied for current
surplus $g^{\delta}(\mathbf{m})\in\mathcal{G}^{\delta} $ if $m_{i}>0$.

\item Control action $\mathbf{E}_{s}$: The manager opts to switch immediately
and the process stops.
\end{itemize}

We denote the space of controls as
\[
\mathcal{E}=\{\mathbf{E}_{s},\left(  \mathbf{E}_{i}\right)  _{i=1,...,n}%
,\mathbf{E}_{0}\}.
\]

Consider $\Pi_{g^{\delta}(\mathbf{m})}^{\delta}\subset\Pi_{g^{\delta
}(\mathbf{m})}$ as the set of all the admissible strategies with initial
surplus $g^{\delta}(\mathbf{m})\in\mathcal{G}^{\delta}$ which can be obtained
by a sequence of control actions in $\mathcal{E}$ at each point of the grid.
Let us describe the strategies $\pi=(\mathbf{L},\overline{\tau})\in
\Pi_{g^{\delta}(\mathbf{m})}^{\delta}$; we take, for any $\omega=(\tau
_{j},\mathbf{U}_{j})_{j\geq1}\in\Omega$, a sequence $\mathbf{s}=(s_{k}%
)_{k=1,...,\tilde{k}}$ with $s_{k}\in\mathcal{E}\ $and $1\leq\tilde{k}%
\leq\infty$, the first control action $s_{1}$ is applied at the point
$g^{\delta}(\mathbf{m})\in\mathcal{G}^{\delta},$ the second control action
$s_{2}$ is applied at the end surplus in $\mathcal{G}^{\delta}$ resulting from
the control action $s_{1},$ and so on. If the length of the sequence
$\mathbf{s}$ is $\tilde{k}<\infty$, then $s_{\tilde{k}}$ should be either
$\mathbf{E}_{s}$ or $\mathbf{E}_{0}$. In the last case, the end surplus
resulting from the final control action $s_{\tilde{k}}$ is outside
$\mathbf{R}_{+}^{n}$ due to the arrival of a jump.

Take $\mathbf{m}^{k}\in\mathbf{N}_{0}^{n}\ $in such a way that $g^{\delta
}(\mathbf{m}^{k})$ is the point of $\mathcal{G}^{\delta}$ in which the control
action $s_{k}$ is applied; let $t_{k}$ be the time in which the control action
$s_{k}$ is chosen; let $\Delta_{k}$ be the time elapsed for the control action
$s_{k}$ and let $\mathbf{y}^{k}\in\mathcal{G}^{\delta}\cup\left(
\mathbf{R}_{+}^{n}\right)  ^{c}$ be the end surplus resulting from the control
action $s_{k}$.

\begin{remark}
Let us describe in a precise way the values of $(\mathbf{m}^{k},\Delta
_{k},t_{k},\mathbf{y}^{k})_{k=1,...,\tilde{k}}$.

\begin{itemize}
\item In the case that $s_{k}=\mathbf{E}_{i}$, then $k<\tilde{k},$ $\Delta
_{k}=0,$ $t_{k+1}=t_{k},$ $\mathbf{m}^{k+1}=\mathbf{m}^{k}-\mathbf{e}_{i}$ and
$\mathbf{y}^{k}=g^{\delta}(\mathbf{m}^{k+1}).$

\item In the case that $s_{k}=\mathbf{E}_{s}$, then
\[
k=\tilde{k},\text{ }t_{k}=\overline{\tau},\Delta_{k}=0\text{ and }%
\mathbf{y}^{k}=g^{\delta}(\mathbf{m}^{k}).
\]

\item In the case that $s_{k}=\mathbf{E}_{0}$, take $j_{k}:=\min\{j:\tau
_{j}>t_{k}\}$ (so $\tau_{j_{k}}$ is the arrival time of the first jump after
$t_{k}$); there are three possibilities:

(a) If $\tau_{j_{k}}>t_{k}+\delta$ , then
\[
k<\tilde{k},\text{ }\Delta_{k}=\delta,\text{ }t_{k+1}=t_{k}+\delta,\text{
}\mathbf{m}^{k+1}=\mathbf{m}^{k}+\mathbf{1}\text{ and }\mathbf{y}%
^{k}=g^{\delta}\left(  \mathbf{m}^{k+1}\right)  .
\]

(b) If $\tau_{j_{k}}\leq t_{k}+\delta$ and $g^{\delta}(\mathbf{m})+\left(
\tau_{j_{k}}-t_{k}\right)  \mathbf{p}-\mathbf{U}_{j}\in\mathbf{R}_{+}^{n}$,
then
\[
k<\tilde{k},~\Delta_{k}=\tau_{j_{k}}-t_{k},~t_{k+1}=\tau_{j_{k}}%
,~\mathbf{m}^{k+1}=\rho^{\delta}(g^{\delta}(\mathbf{m})+\left(  \tau_{j_{k}%
}-t_{k}\right)  \mathbf{p}-\mathbf{U}_{j})~\text{and~}\mathbf{y}^{k}%
=g^{\delta}\left(  \mathbf{m}^{k+1}\right)  .
\]

(c) If $\tau_{j_{k}}\leq t_{k}+\delta$ and $\mathbf{y}^{k}=g^{\delta
}(\mathbf{m})+\left(  \tau_{j_{k}}-t_{k}\right)  \mathbf{p}-\mathbf{U}%
_{j}\notin\mathbf{R}_{+}^{n}$, then
\[
k=\tilde{k},~\Delta_{k}=\tau_{j_{k}}-t_{k}\ ~\text{and}\ ~t_{k}+\Delta
_{k}=\tau_{j_{k}}=\tau^{\mathbf{L}}.
\]

\end{itemize}
\end{remark}

Defining $\Delta\mathbf{L}_{k}$ as the amount of dividends paid by the control
action $s_{k}$, we have%
\[
\Delta\mathbf{L}_{k}=\left\{
\begin{array}
[c]{ll}%
p_{i}\delta\mathbf{e}_{i} & \text{if }s_{k}=\text{$\mathbf{E}_{i}$}\\
\mathbf{c}_{k}-\left\langle \mathbf{c}_{k}\right\rangle ^{\delta} & \text{if
}s_{k}=\text{$\mathbf{E}_{0},$ }\tau_{j_{k}}\in(t_{k},t_{k}+\delta]\text{ and
}\mathbf{c}_{k}\in\mathbf{R}_{+}^{n}\\
0 & \text{otherwise},
\end{array}
\right.
\]
where $j_{k}\ $is defined in the previous remark and
\[
\mathbf{c}_{k}=g^{\delta}(\mathbf{m}^{k})+(\tau_{j_{k}}-t_{k})\mathbf{p}%
-\mathbf{U}_{j}.
\]
Therefore, if the strategy $\pi=(\mathbf{L},\overline{\tau})\in\Pi_{g^{\delta
}(\mathbf{m})}^{\delta}$ then the cumulative dividend payment strategy is
\[
\mathbf{L}_{t}=\sum\nolimits_{k\leq\tilde{k},t_{k}\leq t}\Delta\mathbf{L}_{k},
\]
and the switch time $\overline{\tau}$ is the time in which the control action
$\mathbf{E}_{s}$ is chosen. By construction, if $\pi\in\Pi_{g^{\delta
}(\mathbf{m})}^{\delta}$ then $\mathbf{X}_{t_{k}}^{\mathbf{L}}\in
\mathcal{G}^{\delta}$ for all $k\leq\tilde{k}$ , also the set of times
\[
\{t_{k}:k\leq\tilde{k}\}\subseteq\{\tau_{i}+j\delta:i,j\in\mathbf{N}_{0}\text{
and }j\leq\frac{\tau_{i+1}-\tau_{i}}{\delta}\};
\]
here $\tau_{0}=0$. For the strategy $(\mathbf{L},\overline{\tau})$ to be
admissible, we need to assume the following condition: If the arrival times
and sizes of the claims of two elements in $\Omega$ coincide up to time $t$,
then the corresponding sequences of control actions $\mathbf{s}=(s_{k}%
)_{k=1,...,\tilde{k}}$ must coincide for all $k$ such that $t_{k}\leq t$.

The following lemma states that the sequences $\left(  t_{k}\right)  _{k\geq1}
$ do not have an accumulation point, the proof is in the Appendix.

\begin{lemma}
\label{Lema de tiempo infinito} Given $\pi\in\Pi_{g^{\delta}(\mathbf{m}%
)}^{\delta},$ $\lim_{k\rightarrow\infty}t_{k}=\infty$ a.s. within the subset
$\{\tilde{k}=\infty\}\subset\Omega$.
\end{lemma}

We define the $\mathcal{G}^{\delta}$\textit{-optimal function} $v^{\delta}$ as
the supremum of the value functions of admissible strategies which are
combination of the control actions in $\mathcal{E}$, that is
\begin{equation}
v^{\delta}(\mathbf{m})=\sup\nolimits_{\pi\in\Pi_{g^{\delta}(\mathbf{m}%
)}^{\delta}}V_{\pi}(g^{\delta}(\mathbf{m}))\text{.}\label{vdelta}%
\end{equation}

\subsection{Characterization of the $\mathcal{G}^{\delta}$\textit{-}optimal
function}

In this subsection, we show that the $\mathcal{G}^{\delta}$-optimal function
$v^{\delta}:\mathbf{N}_{0}^{n}\rightarrow\mathbf{R}$ is a solution of a
discrete version of the HJB equation (\ref{HJB}). We also see that $v^{\delta
}$ can be characterized as the smallest supersolution of this discrete HJB
equation. Moreover, we prove that there exists an optimal admissible strategy
for the problem (\ref{vdelta}). This strategy, called the $\mathcal{G}%
^{\delta}$\textit{-optimal strategy}, is stationary in the sense that the
control actions depend only on which point of the grid the current surplus lies.

We now introduce some operators related to the control actions in
$\mathcal{E}$, these operators will be involved in the discrete version of the
HJB equation. Given any family of admissible strategies $\widetilde{\pi}%
=(\pi_{g^{\delta}(\mathbf{m})})_{\mathbf{m}\in\mathbf{N}_{0}^{n}}$ with
$\pi_{g^{\delta}(\mathbf{m})}\in\Pi_{g^{\delta}(\mathbf{m})}^{\delta}$, we
define the value function $\widetilde{w}:\mathbf{N}_{0}^{n}\rightarrow
\mathbf{R}$ of $\widetilde{\pi}$ as
\[
\widetilde{w}(\mathbf{m}):=V_{\pi_{g^{\delta}(\mathbf{m})}}(g^{\delta
}(\mathbf{m})).
\]

Let us consider the admissible strategies with initial surplus $g^{\delta
}(\mathbf{m})\in\mathcal{G}^{\delta}$ which consists on applying first one of
the control actions in $\mathcal{E}$\textbf{,} and afterwards applying the
strategy in the family $\widetilde{\pi}$ corresponding to the end surplus (if
it is possible); the value functions of these strategies are given by%

\begin{equation}%
\begin{array}
[c]{lll}%
T_{0}(\widetilde{w})(\mathbf{m})\bigskip & := & \widetilde{w}(\mathbf{m}%
+\mathbf{1})e^{-(c+\lambda)\delta}+\mathcal{I}^{\delta}(\widetilde
{w})(\mathbf{m})-\int_{0}^{\delta}e^{-(c+\lambda)t}\mathcal{R}(g^{\delta
}(\mathbf{m})+t\mathbf{p})dt\text{ , }%
\end{array}
\label{Definicion TE}%
\end{equation}%
\begin{equation}%
\begin{array}
[c]{ccc}%
T_{i}(\widetilde{w})(\mathbf{m}):=\widetilde{w}(\mathbf{m}-\mathbf{e}%
_{i})+\delta a_{i}p_{i} & \text{and} & T_{s}(\widetilde{w})(\mathbf{m}%
):=f(g^{\delta}(\mathbf{m}))\text{,}%
\end{array}
\label{Definicion Ti TS}%
\end{equation}
depending on which control action in $\mathcal{E}$ is chosen. Here,%
\begin{equation}%
\begin{array}
[c]{l}%
\mathcal{I}^{\delta}(w)(\mathbf{m})\\%
\begin{array}
[c]{ll}%
:= & \int\limits_{0}^{\delta}(%
{\textstyle\int\limits_{\mathbf{0}}^{g^{\delta}(\mathbf{m})+t\mathbf{p}}}
\lambda e^{-(c+\lambda)t}w(\text{$\rho^{\delta}$}(g^{\delta}(\mathbf{m}%
)+t\mathbf{p}-\mathbf{\alpha}))dF(\mathbf{\alpha}))dt\\
& +\int\limits_{0}^{\delta}(%
{\textstyle\int\limits_{\mathbf{0}}^{g^{\delta}(\mathbf{m})+t\mathbf{p}}}
\lambda e^{-(c+\lambda)t}\mathbf{a}\cdot(g^{\delta}(\mathbf{m})+t\mathbf{p}%
-\mathbf{\alpha}-\left\langle g^{\delta}(\mathbf{m})+t\mathbf{p}%
-\mathbf{\alpha}\right\rangle ^{\delta})dF(\mathbf{\alpha}))dt.
\end{array}
\end{array}
\label{Definicion Idelta(W)}%
\end{equation}
We can consider $T_{0}$, $T_{i}$ and $T_{s}$ as operators in the set of
functions $\left\{  w:\mathbf{N}_{0}^{n}\rightarrow\mathbf{R}\right\}  $; we
also define the operator $T$ as%
\begin{equation}
T:=\max\{T_{0},\left(  T_{i}\right)  _{i=1,...,n},T_{s}\}.\label{Definicion T}%
\end{equation}

The following lemma is technical and the proof is in the Appendix.

\begin{lemma}
\label{Ts crecientes} The operators $T_{0}$, $T_{i},$ $T_{s}$ and $T$ are
non-decreasing and $T$ satisfies,
\[
\sup\nolimits_{\mathbf{m}\in\mathbf{N}_{0}^{n}}\left\vert T(w_{1}%
)(\mathbf{m})-T(w_{2})(\mathbf{m})\right\vert \leq\sup\nolimits_{\mathbf{m}%
\in\mathbf{N}_{0}^{n}}\left\vert w_{1}(\mathbf{m})-w_{2}(\mathbf{m}%
)\right\vert .
\]
Moreover, $T_{0}(w)$, $T_{i}(w)$ and $T_{s}(w)$ can be written as a linear
combination of the values of $w(\mathbf{m})$ with $\mathbf{m\in N}_{0}^{n}$
plus a constant.
\end{lemma}

We define the \textit{discrete HJB equation} as%
\begin{equation}
\left(  T(w)-w\right)  (\mathbf{m})=\max\{T_{0}(w)-w,\left(  T_{i}%
(w)-w\right)  _{i=1,...,n},T_{s}(w)-w\}(\mathbf{m})=0\label{Delta HJB}%
\end{equation}
for $\mathbf{m}\in\mathbf{N}_{0}^{n}$. Analogously to Definition
\ref{NuevaDefinicionSubySuper}, we say that a function $\overline
{w}:\mathbf{N}_{0}^{n}\rightarrow\mathbf{R}$ is a \textit{supersolution} of
(\ref{Delta HJB}) if $T(\overline{w})-\overline{w}\leq0$, and a function
$\underline{w}:\mathbf{N}_{0}^{n}\rightarrow\mathbf{R}$ is a
\textit{subsolution} of (\ref{Delta HJB}) if $T(\underline{w})-\underline
{w}\geq0$.

The following results are the discrete versions of Propositions
\ref{Prop V is a viscosity solution}, Lemma
\ref{SupersolucionMayor-ValueFunction}, Theorems \ref{characterization} and
\ref{verification result}. The discrete version of the growth condition
(\ref{gc}) is given by%
\begin{equation}
w(\mathbf{m})e^{\frac{-c}{2n}%
{\textstyle\sum_{i=1}^{n}}
\delta m_{i}}\ \text{is upper bounded in }\mathbf{N}_{0}^{n}\text{.}%
\label{Growth Condition discreta}%
\end{equation}

\begin{proposition}
\label{Propiedades Vdelta}$~$The function $v^{\delta}:\mathbf{N}_{0}%
^{n}\rightarrow\mathbf{R}$ is well defined and it is a solution of
(\ref{Delta HJB}). Moreover, given an initial surplus $g^{\delta}\left(
\mathbf{m}_{0}\right)  \in\mathcal{G}^{\delta}$, there exists a $\mathcal{G}%
^{\delta}$-optimal strategy $\pi_{g^{\delta}\left(  \mathbf{m}_{0}\right)
}^{\delta}\in\Pi_{g^{\delta}\left(  \mathbf{m}_{0}\right)  }^{\delta}$ such
that
\[
v^{\delta}(\mathbf{m}_{0})=V_{\pi_{g^{\delta}\left(  \mathbf{m}_{0}\right)
}^{\delta}}(g^{\delta}\left(  \mathbf{m}_{0}\right)  ).
\]
This $\mathcal{G}^{\delta}$-optimal strategy is\textit{\ stationary} in the
following sense: the control action $s_{k}$ in the sequence $\mathbf{s}%
=(s_{k})_{k=1,...,\tilde{k}}$ depends only on the current surplus $g^{\delta
}\left(  \mathbf{m}^{k}\right)  \in\mathcal{G}^{\delta}$.
\end{proposition}

\textit{Proof}. By definitions (\ref{Definicion V}) and (\ref{vdelta}), we
have
\[
f(g^{\delta}(\mathbf{m}))\leq v^{\delta}(\mathbf{m})\leq V(g^{\delta
}(\mathbf{m})),
\]
so $v^{\delta}$ is well defined.

Let us prove that $v^{\delta}=T(v^{\delta})$. Take a sequence $(p_{l}%
)_{l\geq1}$ of families of strategies $p_{l}=(\pi_{g^{\delta}(\mathbf{m})}%
^{l})_{\mathbf{m}\in\mathbf{N}_{0}^{n}}$ with $\pi_{g^{\delta}(\mathbf{m}%
)}^{l}\in\Pi_{g^{\delta}(\mathbf{m})}^{\delta}$ such that
\[
v^{\delta}(\mathbf{m})-V_{\pi_{g^{\delta}(\mathbf{m})}^{l}}(g^{\delta
}(\mathbf{m}))\leq\frac{1}{l}%
\]
for all $\mathbf{m}\in\mathbf{N}_{0}^{n}$. Define $w_{l}:\mathbf{N}_{0}%
^{n}\rightarrow\mathbf{R}$ as $w_{l}(\mathbf{m})=V_{\pi_{g^{\delta}%
(\mathbf{m})}^{l}}(g^{\delta}(\mathbf{m}))$, by Lemma \ref{Ts crecientes}, we
have that
\[
T(v^{\delta})(\mathbf{m})=\lim_{l\rightarrow\infty}T(w_{l})(\mathbf{m})\leq
v^{\delta}(\mathbf{m})\text{.}%
\]
On the other hand, since $\pi_{g^{\delta}(\mathbf{m})}^{l}$ can be obtained by
a sequence of control actions $\mathbf{s}=(s_{k})_{k=1,...,\tilde{k}}$ and at
any point $g^{\delta}(\mathbf{m})$ of the grid all the value functions of
strategies in $\Pi_{g^{\delta}(\mathbf{m})}^{\delta}$ are bellow $v^{\delta
}(\mathbf{m})$, we have by definition of $T$ given in (\ref{Definicion T}),
that $w_{l}(\mathbf{m})\leq T(v^{\delta})(\mathbf{m})$. So taking the limit as
$l\rightarrow\infty,$ we obtain that
\[
v^{\delta}(\mathbf{m})\leq T(v^{\delta})(\mathbf{m}).
\]

Finally, since $v^{\delta}=T(v^{\delta})$, we can define for any
$\mathbf{m}\in\mathbf{N}_{0}^{n}$, a control action $S(\mathbf{m}%
)\in\mathcal{E}$ in the following way:

\begin{itemize}
\item If $T_{s}(v^{\delta})(\mathbf{m})=v^{\delta}(\mathbf{m})$, take
$S(\mathbf{m})=\mathbf{E}_{s}$.

\item If $T_{0}(v^{\delta})(\mathbf{m})=v^{\delta}(\mathbf{m})$, take
$S(\mathbf{m})=\mathbf{E}_{0}$.

\item and if $T_{i}(v^{\delta})(\mathbf{m})=v^{\delta}(\mathbf{m})$ for some
$i=1,...,n$, take $S(\mathbf{m})=\mathbf{E}_{i}$\textbf{.}
\end{itemize}

Given an initial surplus $g^{\delta}(\mathbf{m}_{0})\in\mathcal{G}^{\delta},$
the $\mathcal{G}^{\delta}$-optimal\textit{\ }strategy\textit{\ }%
$\pi_{g^{\delta}(\mathbf{m}_{0})}^{\delta}\in\Pi_{g^{\delta}(\mathbf{m}_{0}%
)}^{\delta}$ is defined inductively as follows: $s_{1}=S(\mathbf{m}_{0}) $;
assuming that $s_{1},s_{2},..,s_{k-1}$ are defined and the process does not
stop at step $k-1$, we define $s_{k}=S(\mathbf{m}_{0}^{k})$ where $g^{\delta
}\left(  \mathbf{m}_{0}^{k}\right)  \in$ $\mathcal{G}^{\delta}$ is the end
surplus of $s_{k-1}$. $\blacksquare$

Analogously to Remark \ref{Muchas soluciones viscosas}, we cannot expect in
general to have uniqueness of viscosity solutions of the discrete HJB equation
(\ref{Delta HJB}). For instance, in the two dimensional case with independent
companies, the switch function $f$ given in (\ref{Merger 2x2}) and the
function $\upsilon$ defined in (\ref{Nu the dos companias independientes}), we
have that
\[
w(\mathbf{m}):=%
{\textstyle\sum_{i=1}^{n}}
p_{i}m_{i}\delta+k
\]
is a solution of (\ref{Delta HJB}) for $k$ large enough. The following lemma
is the discrete version of Lemma \ref{SupersolucionMayor-ValueFunction}, the
proof is in the Appendix.

\begin{lemma}
\label{Menor Supersolucion Discreta} Given any $\pi=(\mathbf{L},\overline
{\tau})\in$ $\Pi_{g^{\delta}(\mathbf{m})}^{\delta}$ and any supersolution $w:$
$\mathbf{N}_{0}^{n}\rightarrow\mathbf{R}$ of (\ref{Delta HJB}) with growth
condition (\ref{Growth Condition discreta}), we have that $V_{\pi}(g^{\delta
}(\mathbf{m}))\leq w(\mathbf{m})$.
\end{lemma}

From Lemma \ref{Menor Supersolucion Discreta}, we obtain the following theorems.

\begin{theorem}
The $\mathcal{G}^{\delta}$-optimal value function $v^{\delta}$ $:\mathbf{N}%
_{0}^{n}\rightarrow\mathbf{R}$ can be characterized as the smallest
supersolution of the discrete HJB equation (\ref{Delta HJB}) with growth
condition (\ref{Growth Condition discreta}).
\end{theorem}

\begin{theorem}
\label{TeoremaVerificacionDiscreto} If the function $w:$ $\mathbf{N}_{0}%
^{n}\rightarrow\mathbf{R}$ with growth condition
(\ref{Growth Condition discreta})$\ $is a supersolution of (\ref{Delta HJB}),
and also satisfies that for any $\mathbf{m}\in\mathbf{N}_{0}^{n},$
$w(\mathbf{m})$ is either $V_{\pi}(g^{\delta}(\mathbf{m}))$ with $\pi\in
\Pi_{g^{\delta}(\mathbf{m})}^{\delta}$ or $\lim_{l\rightarrow\infty}V_{\pi
_{l}}(g^{\delta}(\mathbf{m}))$ with $\pi_{l}\in\Pi_{g^{\delta}(\mathbf{m}%
)}^{\delta}$ for any $l\geq1$, then $w=v^{\delta}$.
\end{theorem}

\subsection{Construction of the $\mathcal{G}^{\delta}$\textit{-}optimal
strategy and the $\mathcal{G}^{\delta}$\textit{-}optimal function}

In this subsection we construct recursively the $\mathcal{G}^{\delta}$-optimal
strategy and the corresponding $\mathcal{G}^{\delta}$-optimal function.

Since $T$ defined in (\ref{Definicion T}) is not a contraction operator,
$v^{\delta}$ can not be obtained numerically as a fixed point; so we construct
value functions $v_{l}^{\delta}$ of strategies in $\Pi_{g^{\delta}%
(\mathbf{m})}^{\delta}$ which can be calculated explicitly by
(\ref{Definicion TE}), (\ref{Definicion Ti TS}) and (\ref{Definicion T}) such
that $v_{l}^{\delta}$ $\nearrow$ $v^{\delta}$ as $l\rightarrow\infty$.

Let us define iteratively the families of strategies $\widetilde{\pi}_{l}%
=(\pi_{g^{\delta}(\mathbf{m})}^{l})_{\mathbf{m}\in\mathbf{N}_{0}^{n}}$ for
each $l\geq1$ in the following way:

(1) We start with the family of strategies $\widetilde{\pi}_{1}=(\pi
_{g^{\delta}(\mathbf{m})}^{1})_{\mathbf{m}\in\mathbf{N}_{0}^{n}}$ where
$\pi_{g^{\delta}(\mathbf{m})}^{1}\in\Pi_{g^{\delta}(\mathbf{m})}^{\delta}$
consists on switching immediately; the value of this strategy is
\[
v_{1}^{\delta}(\mathbf{m}):=f(g^{\delta}(\mathbf{m})).
\]

(2) Given the family of strategies $\widetilde{\pi}_{l}=(\pi_{g^{\delta
}(\mathbf{m})}^{l})_{\mathbf{m}\in\mathbf{N}_{0}^{n}}$ with $\pi_{g^{\delta
}(\mathbf{m})}^{l}\in\Pi_{g^{\delta}(\mathbf{m})}^{\delta}$, we define the
family $\widetilde{\pi}_{l+1}=(\pi_{g^{\delta}(\mathbf{m})}^{l+1}%
)_{\mathbf{m}\in\mathbf{N}_{0}^{n}}$ as follows: We choose for any
$\mathbf{m}\in\mathbf{N}_{0}^{n}$, the best strategy $\pi_{g^{\delta
}(\mathbf{m})}^{l+1}\in\Pi_{g^{\delta}(\mathbf{m})}^{\delta}$ among the ones
which initially follows one of control actions in $\mathcal{E}$ and then
continues with the corresponding strategy in the family $\widetilde{\pi}_{l}$.
The value of this new strategy is given by%
\begin{equation}
v_{l+1}^{\delta}(\mathbf{m}):=T(v_{l}^{\delta})(\mathbf{m})=T^{l}%
(v_{1}^{\delta})(\mathbf{m})\text{ for }\mathbf{m}\in\mathbf{N}_{0}%
^{n}.\label{Definicionvk}%
\end{equation}

\begin{remark}
$v_{l}^{\delta}$ can be thought as the maximum of the value function of
strategies $\pi\in\Pi_{g^{\delta}(\mathbf{m})}^{\delta}$ where the length
$\tilde{k}$ of the corresponding sequence $\mathbf{s}\ $is upper bounded by
$l$ and $s_{l}=\mathbf{E}_{s}\mathbf{\ }$in the case that $\tilde{k}=l.$
\end{remark}

In the next proposition we use Theorem \ref{TeoremaVerificacionDiscreto} to
see that the limit of $v_{l}^{\delta}$ is indeed $v^{\delta}$.

\begin{proposition}
We have that $v_{l+1}^{\delta}\geq v_{l}^{\delta}$ for $l\geq1$ and that
$\lim_{l\rightarrow\infty}v_{l}^{\delta}=v^{\delta}.$
\end{proposition}

\textit{Proof}. Take $\mathbf{m}\in\mathbf{N}_{0}^{n}$, it is straightforward
to see by (\ref{Definicionvk}) that $v_{2}^{\delta}(\mathbf{m})\geq
v_{1}^{\delta}(\mathbf{m})$; on the other hand, the operator $T$ is
non-decreasing, so we obtain that $v_{l+1}^{\delta}\geq v_{l}^{\delta}$ for
$l\geq1$. Then, there exists $w_{0}:\mathbf{N}_{0}^{n}$ $\rightarrow
\mathbf{R}$ such that
\[
w_{0}(\mathbf{m}):=\lim\nolimits_{l\rightarrow\infty}v_{l}^{\delta}%
(\mathbf{m})\leq V(g^{\delta}(\mathbf{m})).
\]
Note that all the functions $v_{l}^{\delta}$ are subsolutions (\ref{Delta HJB}%
) and that $w_{0}$ is a solution of (\ref{Delta HJB}) because $T(w_{0})=w_{0}%
$. Since $w_{0}$ satisfies the growth condition
(\ref{Growth Condition discreta}), $w_{0}$ coincides with the value function
$v^{\delta}$ by Theorem \ref{TeoremaVerificacionDiscreto}. $\blacksquare$

\subsection{Definition of the value function $V^{\delta}$}

In this subsection we define, using the $\mathcal{G}^{\delta}$-optimal
functions and strategies, a family of admissible strategies for any point in
$\mathbf{R}_{+}^{n}$ and the corresponding value function $V^{\delta}.$

\begin{description}
\item
\begin{definition}
\label{Definicion de Vdelta en la grilla} We use the $\mathcal{G}^{\delta}%
$-optimal function $v^{\delta}:\mathbf{N}_{0}^{n}\rightarrow\mathbf{R}$ to
define a function $V^{\delta}:\mathcal{G}^{\delta}\rightarrow\mathbf{R}$ as
\[
V^{\delta}(g^{\delta}(\mathbf{m})):=v^{\delta}(\mathbf{m})
\]
for $\mathbf{m}\in\mathbf{N}_{0}^{n}$. Note that $V^{\delta}(g^{\delta
}(\mathbf{m}))$ is the value of the $\mathcal{G}^{\delta}$-optimal admissible
strategy $\pi_{g^{\delta}(\mathbf{m})}^{\delta}\in\Pi_{g^{\delta}(\mathbf{m}%
)}^{\delta}$.
\end{definition}
\end{description}

We construct now a family of strategies $\widetilde{\pi}^{\delta}=\left(
\pi_{\mathbf{x}}\right)  _{\mathbf{x}\in\mathbf{R}_{+}^{n}}$, where
$\pi_{\mathbf{x}}\in\Pi_{\mathbf{x}}$, such that the corresponding value
function $V^{\delta}(\mathbf{x})=V_{\pi_{\mathbf{x}}}(\mathbf{x})$ extends to
$\mathbf{R}_{+}^{n}$ the function defined in Definition
\ref{Definicion de Vdelta en la grilla}. Take the strategy $\pi_{\mathbf{x}%
}\in\Pi_{\mathbf{x}}$ which pays immediately $\mathbf{x}-\left\langle
\mathbf{x}\right\rangle ^{\delta}$ as dividends and then follows the
$\mathcal{G}^{\delta}$-optimal strategy $\pi_{\left\langle \mathbf{x}%
\right\rangle ^{\delta}}^{\delta}\in\Pi_{\left\langle \mathbf{x}\right\rangle
^{\delta}}^{\delta}$. We obtain that $V^{\delta}:\mathbf{R}_{+}^{n}%
\rightarrow\mathbf{R}$ is given by
\begin{equation}
V^{\delta}(\mathbf{x}):=V^{\delta}(\left\langle \mathbf{x}\right\rangle
^{\delta})+\mathbf{a}\cdot(\mathbf{x}-\left\langle \mathbf{x}\right\rangle
^{\delta}).\label{Definicion Vdelta}%
\end{equation}

\section{Convergence of the Discrete Scheme}

In this section we show the locally uniformly convergence of the discrete
scheme defined in the previous section by taking a suitable sequence of
embedded grids.

In the next technical lemma, we show that the functions $v^{\delta}$ satisfy a
$\delta$-locally Lipschitz condition and a relation between $v^{2\delta}$ and
$v^{\delta}$ which gives a monotonicity condition on the embedded grids; the
proof is in the Appendix.

\begin{lemma}
\label{v_delta es lipschitz y monotona en el reticulado}The functions
$v^{\delta}$ defined in (\ref{vdelta}) satisfy:

(1) $v^{\delta}(\mathbf{m}+\mathbf{e}_{i})-v^{\delta}(\mathbf{m})\geq
a_{i}p_{i}\delta$ and $v^{\delta}(\mathbf{m}+\mathbf{1})-v^{\delta}%
(\mathbf{m})\leq v^{\delta}(\mathbf{m})(e^{(c+\lambda)\delta}-1);$

(2) $\Pi_{g^{2\delta}(\mathbf{m})}^{2\delta}\subset\Pi_{2g^{\delta}%
(\mathbf{m})}^{\delta}\subset\Pi_{2g^{\delta}(\mathbf{m})}$ and so
$v^{2\delta}(\mathbf{m})\leq v^{\delta}(2\mathbf{m})$.
\end{lemma}

Let us take $\delta_{k}:=\delta/2^{k}$ for $k\geq0$. In the remainder of the
section we will prove that $V^{\delta_{k}}$ $\nearrow V$ locally uniformly as
$k$ goes to infinity. Consider the dense set in $\mathbf{R}_{+}^{n}$,
$\mathcal{G}:=\bigcup\nolimits_{k\geq0}\mathcal{G}^{\delta_{k}}$. Note that
$\mathcal{G}^{\delta_{k}}\subset\mathcal{G}^{\delta_{k+1}}$, so by Lemma
\ref{v_delta es lipschitz y monotona en el reticulado}-(2),
\[
V^{\delta_{k}}\leq V^{\delta_{k+1}}\leq V;
\]
then we can define the function $\overline{V}:$ $\mathbf{R}_{+}^{n}%
\rightarrow\mathbf{R}$ as%

\begin{equation}
\overline{V}(\mathbf{x}):=\lim\nolimits_{k\rightarrow\infty}V^{\delta_{k}%
}(\mathbf{x}).\label{ubarra como limite}%
\end{equation}

\begin{remark}
\label{Limit value function u barra} We will prove that $\overline{V}$ is the
optimal value function. In order to do that, we will show that $\overline{V}$
is a viscosity supersolution of (\ref{HJB}). It is straightforward to see that
$\overline{V}(\mathbf{x})$ is a limit of value functions of admissible
strategies in $\Pi_{\mathbf{x}}$ for all $\mathbf{x}\in\mathbf{R}_{+}^{n}$ so
the result will follow from Theorem \ref{verification result}. Since there is
no uniqueness of solution of the HJB equation, it is essential to show that
this function is a limit of value functions of admissible strategies.
\end{remark}

In the next lemma, we find a bound on the variation of $V^{\delta_{k}}$ and we
show that $\overline{V}$ is locally Lipschitz in $\mathbf{R}_{+}^{n}$ and so
it is absolutely continuous; the proof is in the Appendix.

\begin{lemma}
\label{Lipschitz Inequality u barra} We have for each $\mathbf{y}%
\geq\mathbf{x}$ in $\mathbf{R}_{+}^{n}$ that
\[
\left\vert V^{\delta_{k}}(\mathbf{y})-V^{\delta_{k}}(\mathbf{x})\right\vert
\leq\left\Vert \left\langle \mathbf{y}\right\rangle ^{\delta_{k}}-\left\langle
\mathbf{x}\right\rangle ^{\delta_{k}}\right\Vert _{1}\frac{2}{\hat{p}%
}V^{\delta_{k}}(\left\langle \mathbf{x}\vee\mathbf{y}\right\rangle
^{\delta_{k}})(\frac{e^{(c+\lambda)\delta_{k}}-1}{\delta_{k}})+2\delta
_{k}\mathbf{a}\cdot\mathbf{p}\text{,}%
\]
and also%
\[
\mathbf{a}\cdot\left(  \mathbf{y}-\mathbf{x}\right)  \leq\overline
{V}(\mathbf{y})-\overline{V}(\mathbf{x})\leq\overline{V}(\mathbf{y}%
)\frac{2(c+\lambda)}{\hat{p}}\left\Vert \mathbf{y}-\mathbf{x}\right\Vert
_{1}\text{,}%
\]
where $\hat{p}:=\min_{i=1,..,n}p_{i}$.
\end{lemma}

In the next two propositions we address the convergence of $V^{\delta_{k}}$ to
$\overline{V}$ and we prove that $\overline{V}$ coincides with $V$.

\begin{proposition}
\label{Limite u barra uniforme}For any $\delta>0$, $V^{\delta_{k}}$
$\nearrow\overline{V}\ $ locally uniformly as $k$ goes to infinity.
\end{proposition}

\textit{Proof}. Consider a compact set $K$ in $\mathbf{R}_{+}^{n}$,
$\mathbf{x}^{1}\in K$ and $\varepsilon>0$. Let us take\ an upper bound
$\mathbf{z}\in\mathbf{R}_{+}^{n}$ of $K$. We show first that there exists
$k_{0}$ large enough and $\eta>0$ small enough such that if $\left\Vert
\mathbf{x}-\mathbf{x}^{1}\right\Vert _{1}<\eta$ and $k\geq$ $k_{1}$, then%
\begin{equation}
\overline{V}(\mathbf{x})-V^{\delta_{k}}(\mathbf{x})<\varepsilon
.\label{diferencia en Bolas}%
\end{equation}
Indeed, by pointwise convergence at $\mathbf{x}^{1}$, there exists $k_{1}$
such that
\begin{equation}
\overline{V}(\mathbf{x}^{1})-V^{\delta_{k}}(\mathbf{x}^{1})<\varepsilon
/3~\text{for }k\geq k_{1}.\label{desig 1}%
\end{equation}
By Lemma \ref{Lipschitz Inequality u barra}, there exists $\eta_{1}$ such that
if $\left\Vert \mathbf{x}-\mathbf{x}^{1}\right\Vert _{1}<\eta_{1}$, then
\begin{equation}
\left\vert \overline{V}(\mathbf{x})-\overline{V}(\mathbf{x}^{1})\right\vert
<\varepsilon/3.\label{desig 2}%
\end{equation}
Also, from Lemma \ref{Lipschitz Inequality u barra}, there exists $\eta_{2}$
and $k_{2}$ such that if $\left\Vert \mathbf{x}-\mathbf{x}^{1}\right\Vert
_{1}<\eta_{1}$, then%
\begin{equation}
\left\vert V^{\delta_{k}}(\mathbf{x})-V^{\delta_{k}}(\mathbf{x}^{1}%
)\right\vert \leq\left\Vert g^{\delta_{k}}\left(  \rho^{\delta_{k}}%
(\mathbf{x})-\rho^{\delta_{k}}(\mathbf{x}^{1})\right)  \right\Vert
_{1}\overline{V}(\mathbf{z})2e^{(c+\lambda)}/\hat{p}+2\delta_{k}%
\mathbf{a}\cdot\mathbf{p}<\varepsilon/3\label{desig 3}%
\end{equation}
for $k\geq k_{2}$. Therefore, taking $\eta:=\eta_{1}\wedge\eta_{2}$, for
$k\geq k_{0}:=k_{1}\vee k_{2}$, we obtain (\ref{diferencia en Bolas}) from
(\ref{desig 1}), (\ref{desig 2}) and (\ref{desig 3}).

Finally, we conclude the result taking a finite covering of the compact set
$K$. $\blacksquare$

\begin{proposition}
\label{vbarra es supersolucion} The function $\overline{V}$ defined in
(\ref{ubarra como limite}) is the optimal value function $V$.
\end{proposition}

\textit{Proof}. By Remark \ref{Limit value function u barra}, it is enough to
prove that $\overline{V}$ is a viscosity supersolution of (\ref{HJB}) in the
interior of $\mathbf{R}_{+}^{n}$. Take $\mathbf{x}^{0}$ in the interior of
$\mathbf{R}_{+}^{n}$ and a differentiable test function $\varphi
:\mathbf{R}_{+}^{n}\rightarrow\mathbf{R}$ for viscosity supersolution of
(\ref{HJB}) at $\mathbf{x}^{0}$, that is
\begin{equation}
\overline{V}(\mathbf{x})\geq\varphi(\mathbf{x})\text{ and }\overline
{V}(\mathbf{x}^{0})=\varphi(\mathbf{x}^{0}).\label{Comparacion}%
\end{equation}
Since $\mathcal{G}$ is a dense set in $\mathbf{R}_{+}^{n}$, we obtain by the
continuity assumptions on the function $f$ given in Section
\ref{Seccion Modelistica} and (\ref{ubarra como limite}) that $f\leq
\overline{V}$ in $\mathbf{R}_{+}^{n}$, so $f(\mathbf{x}^{0})-\varphi
(\mathbf{x}^{0})\leq0$. By Proposition \ref{Lipschitz Inequality u barra},
\[
\overline{V}(\mathbf{y})-\overline{V}(\mathbf{x})\geq\mathbf{a}\cdot\left(
\mathbf{y}-\mathbf{x}\right)
\]
for all $\mathbf{y}\geq\mathbf{x}$, so it holds that $\mathbf{a}-\nabla
\varphi(\mathbf{x}^{0})\leq\mathbf{0}$. In order to prove that $\mathcal{L}%
(\varphi)(\mathbf{x}^{0})\leq0,$ consider now for $\eta>0$ small enough,
\[
\varphi_{\eta}(\mathbf{x})=\varphi(\mathbf{x})-\eta\left(  \mathbf{x}%
-\mathbf{x}^{0}\right)  \mathbf{\cdot}(\mathbf{x}-\mathbf{x}^{0}).
\]
Given $k\geq0$, the set $\mathcal{G}^{\delta_{k}}\cap\lbrack\mathbf{0}%
,\mathbf{x}^{0}+\mathbf{1}]$ is finite, so we can define
\begin{equation}
a_{k}^{\eta}:=\min\nolimits_{\mathcal{G}^{\delta_{k}}\cap\lbrack
\mathbf{0},\mathbf{x}^{0}+\mathbf{1}]}\{V^{\delta_{k}}(\mathbf{x}%
)-\varphi_{\eta}(\mathbf{x})\}.\label{minimo V_delta-Fi}%
\end{equation}
Since $V^{\delta_{k}}\leq\overline{V}$, we have from (\ref{Comparacion}), that
$a_{k}^{\eta}\leq0$. Taking%
\[
0\leq b_{k}:=\max\nolimits_{\mathcal{G}^{\delta_{k}}\cap\lbrack\mathbf{0}%
,\mathbf{x}^{0}+\mathbf{1}]}\left(  \overline{V}-V^{\delta_{k}}\right)  ,
\]
by Proposition \ref{Limite u barra uniforme}, $b_{k}\rightarrow0$ as
$k\rightarrow\infty$. \ For all $\mathbf{x}\in\mathcal{G}^{\delta_{k}}%
\cap\lbrack\mathbf{0},\mathbf{x}^{0}+\mathbf{1}]$ we get from
(\ref{Comparacion}),%

\[%
\begin{array}
[c]{lll}%
V^{\delta_{k}}(\mathbf{x})-\varphi_{\eta}(\mathbf{x}) & = & V^{\delta_{k}%
}(\mathbf{x})-\overline{V}(\mathbf{x})+\overline{V}(\mathbf{x})-\varphi
(\mathbf{x})+\eta\left(  \mathbf{x}-\mathbf{x}^{0}\right)  \mathbf{\cdot
}(\mathbf{x}-\mathbf{x}^{0})\\
& \geq & -b_{k}+\eta\left(  \mathbf{x}-\mathbf{x}^{0}\right)  \mathbf{\cdot
}(\mathbf{x}-\mathbf{x}^{0}).
\end{array}
\]
Then, the minimum argument in (\ref{minimo V_delta-Fi}) is attained at
$\mathbf{x}^{k}\in\mathcal{G}^{\delta_{k}}$ such that
\[
\left(  \mathbf{x}^{k}-\mathbf{x}^{0}\right)  \mathbf{\cdot}(\mathbf{x}%
^{k}-\mathbf{x}^{0})\leq b_{k}/\eta.
\]
Then, we have $\mathbf{x}^{k}\rightarrow\mathbf{x}^{0}$ and $-a_{k}^{\eta
}\rightarrow0$ as $k$ goes to infinity. So
\[
V^{\delta_{k}}(\mathbf{x})\geq\varphi_{\eta}(\mathbf{x})-a_{k}^{\eta}\text{
for }\mathbf{x}\in\mathcal{G}^{\delta_{k}}\cap\lbrack0,\mathbf{x}%
^{0}+\mathbf{1}]\text{ and }V^{\delta_{k}}(\mathbf{x}^{k})=\varphi_{\eta
}(\mathbf{x}^{k})-a_{k}^{\eta}.
\]
Since
\[
T_{0}(v^{\delta_{k}})\left(  \left[  \tfrac{x_{1}^{k}}{\delta_{k}p_{1}%
}\right]  ,...,\left[  \tfrac{x_{n}^{k}}{\delta_{k}p_{n}}\right]  \right)
-v^{\delta_{k}}\left(  \left[  \tfrac{x_{1}^{k}}{\delta_{k}p_{1}}\right]
,...,\left[  \tfrac{x_{n}^{k}}{\delta_{k}p_{n}}\right]  \right)  \leq0,
\]
we obtain
\[%
\begin{array}
[c]{lll}%
0 & \geq & e^{-(c+\lambda)\delta_{k}}\left(  V^{\delta_{k}}(\mathbf{x}%
^{k}+\delta_{k}\mathbf{p})\right) \\
&  & +\int_{0}^{\delta_{k}}\lambda e^{-(c+\lambda)t}(%
{\textstyle\int\nolimits_{\mathbf{0}\leq\mathbf{\alpha}\leq\mathbf{x}%
^{k}+t\mathbf{p}}}
V^{\delta_{k}}(\mathbf{x}^{k}+t\mathbf{p}-\mathbf{\alpha})dF(\mathbf{\alpha
}))dt\\
&  & -\int_{0}^{\delta_{k}}e^{-(c+\lambda)t}\mathcal{R}(\mathbf{x}%
^{k}+t\mathbf{p})dt-V^{\delta_{k}}(\mathbf{x}^{k})\\
& \geq & e^{-(c+\lambda)\delta_{k}}\left(  \varphi_{\eta}(\mathbf{x}%
^{k}+\delta_{k}\mathbf{p})-\varphi_{\eta}(\mathbf{x}^{k})\right) \\
&  & -\left(  \varphi_{\eta}(\mathbf{x}^{k})-a_{k}^{\eta}\right)
(1-e^{-(c+\lambda)\delta_{k}})\\
&  & +\int_{0}^{\delta_{k}}\lambda e^{-(c+\lambda)t}(%
{\textstyle\int\nolimits_{\mathbf{0}\leq\mathbf{\alpha}\leq\mathbf{x}%
^{k}+t\mathbf{p}}}
\left(  \varphi_{\eta}(\rho^{\delta_{k}}(\mathbf{x}^{k}+t\mathbf{p}%
-\mathbf{\alpha})-a_{k}^{\eta}\right)  dF(\mathbf{\alpha}))dt\\
&  & +\int_{0}^{\delta_{k}}\lambda e^{-(c+\lambda)t}(%
{\textstyle\int\nolimits_{\mathbf{0}\leq\mathbf{\alpha}\leq\mathbf{x}%
^{k}+t\mathbf{p}}}
\left(  \mathbf{a}\cdot\left(  \mathbf{x}^{k}+t\mathbf{p}-\mathbf{\alpha
}-\left\langle \mathbf{x}^{k}+t\mathbf{p}-\mathbf{\alpha}\right\rangle
^{\delta_{k}}\right)  \right)  dF(\mathbf{\alpha}))dt\\
&  & -\int_{0}^{\delta_{k}}e^{-(c+\lambda)t}\mathcal{R}(\mathbf{x}%
^{k}+t\mathbf{p})dt\text{.}%
\end{array}
\]
Dividing by $\delta_{k}$, taking $k\ $to infinity and using the continuity of
$\mathcal{R}$, we get $\mathcal{L}(\varphi_{\eta})(\mathbf{x}^{0})\leq0 $.
Finally, since $\nabla\varphi_{\eta}(\mathbf{x}^{0})=\nabla\varphi
(\mathbf{x}^{0})$ and $\varphi_{\eta}\nearrow\varphi$ as $\eta\searrow0$, we
obtain that $\mathcal{L}(\varphi)(\mathbf{x}^{0})\leq0$ and the result
follows. $\blacksquare$

From Propositions \ref{Limite u barra uniforme} and
\ref{vbarra es supersolucion}, we conclude the main result of the paper.

\begin{theorem}
\label{Main Theorem} For any $\delta>0$, the functions $V^{\delta_{k}}$
$\nearrow\overline{V}=V$ locally uniformly as $k$ goes to infinity.
\end{theorem}

\section{Optimal merger time}

Let us assume that the uncontrolled bivariate surplus $\mathbf{X}_{t}$ of two
insurance companies with the same shareholders follows the process
(\ref{UncontrolledSurplusOriginal}). Both branches pay dividends up to the
time of their respective ruin $\tau^{L_{i}}$ with $i=1,2$, but the
shareholders has the possibility of \textit{merging} the two branches at any
time $\overline{\tau}$ prior to $\tau^{\mathbf{L}}=\tau^{L_{1}}\wedge
\tau^{L_{2}}$ (as defined in (\ref{Definicion Tau L})); at this time the
branches put together all their surplus, pay the claims of both branches and
pay dividends until the joined surplus becomes negative, see e.g. Gerber and
Shiu \cite{GS Merger}. The aim is to find both the dividend payment policy and
the merging time which maximize the expected sum of all the discounted
dividends paid to the shareholders. This problem corresponds to
(\ref{Definicion V}) where $n=2$, $\mathbf{a}=(1,1)$, $A$ is the $2\times2$
identity matrix, the function $\upsilon$ is defined as in
(\ref{Nu the dos companias independientes}) and the the switch-value function
$f$ is defined as in (\ref{Merger 2x2}). In the numerical examples, we
consider
\[
F(x_{1},x_{2})=\mathbb{P}(\alpha_{1}\leq x_{1},\alpha_{2}\leq x_{2}%
)=\frac{\lambda_{1}}{\lambda}(1-e^{-d_{1}x_{1}})+\frac{\lambda_{2}}{\lambda
}(1-e^{-d_{2}x_{2}})
\]
with $d_{1}=3$ and $d_{2}=3.5$. Note that the above formula for $F$
corresponds to the case in which the surplus processes of the two branches are
independent, as we pointed out in (\ref{independent_Surpluses}); so the
function
\[
\mathcal{R(}x_{1},x_{2})=\frac{\lambda_{1}}{\lambda}V_{2}(x_{2})e^{-d_{1}%
x_{1}}+\frac{\lambda_{2}}{\lambda}V_{1}(x_{1})e^{-d_{2}x_{2}}%
\]
is continuous in $\mathbf{R}_{+}^{2}$. The parameters of the merger company
(that is a one dimension problem) are $\lambda_{M}=\lambda_{1}+\lambda_{2} $,
$p_{M}=p_{1}+p_{2}$ and $F_{M}(x)=F(x,x)$.

In the first example, we consider $\lambda_{1}=2.4$, $\lambda_{2}=2$,
$\lambda=\lambda_{1}+\lambda_{2}$, $p_{1}=1.08$, $p_{2}=0.674$, $c=0.11$,
$\delta=1/60$ and $c_{M}=0$. In Figure 1, we show the $\mathcal{G}^{\delta}%
$-optimal strategy: the merger region is in black, the non-action region in
white, the dividend payment region for the first company region in dark grey
and the dividend payment region for the second company in light grey. Note
that the non-action region has two connected components; in the one on the
top, the optimal strategy is to withhold dividend payments in order to reach
the merger region, and in the white rectangle on the bottom the optimal
strategy corresponds to the non-action region of the stand-alone problem (in
which the companies never merge). This figure suggests that, as $\delta
\rightarrow0$, the optimal local control in the boundary between the
non-action rectangle and the dividend payment region for the second company
region (light grey), should be that the second company pay the incoming
premium as dividends while the first company pays no dividends, so the
bivariate control surplus stays on the top boundary $x_{2}=0.33$ of the
rectangle and moves rightward at constant speed $p_{1}$ to the point
$(0.33,1.42)$, which corresponds to the righ-top corner of the rectangle
(until the arrival of the next claim). Analogously, the optimal strategy in
the right boundary $x_{1}=1.42$ of the non action rectangle should be that the
first company pay the incoming premium as dividends while the second company
pay no dividends, in this case the bivariate control surplus stays on the
right boundary of the rectangle and moves upward at constant speed $p_{2}$ to
the righ-top corner (until the arrival of the next claim). At this corner,
both companies pay their incoming premium as dividends and the surplus process
remains constant (until the arrival of the next claim). It is more difficult
to guess the optimal local control (as $\delta\rightarrow0)$ in the boundary
between the upper connected component of the non-action region and the
dividend payment region for the second company region (light grey). Our
conjecture, assuming some regularity on this boundary, is the following: In
the upper part of this boundary (up to the furthest point to the right), the
second company should pay dividends with some rate in such a way that the
bivariate control surplus stays in this part of the boundary (moving
downwards), and in the lower part of this boundary, the second company should
pay a lump sum in such a way that the bivariant surplus reaches the line
$x_{2}=0.33$.

In the second example, we consider $\lambda_{1}=2.44$, $\lambda_{2}=2.22$,
$\lambda=\lambda_{1}+\lambda_{2}$, $p_{1}=1.100$, $p_{2}=0.825$, $c=0.1$,
$\delta=1/50$ and $c_{M}=0.364$. In Figure 2, we show the $\mathcal{G}%
^{\delta}$-optimal strategy; the regions are described with the same colors as
before. This figure suggests that, as $\delta\rightarrow0$, the optimal local
control in the boundary between the non-action region (white) and the dividend
payment region for the second company (light grey region), would be (assuming
some regularity on the boundary) that the second company pay dividends with
some rate in such a way that the bivariate control surplus stays in the
boundary: this control surplus would move downward until the bivariate surplus
reach the point $(1.61,1.06)$ in which the light grey, the dark grey and the
white regions meet. At this point, both companies should pay the incoming
premiums as dividends and the bivariate surplus process remains constant until
the arrival of the next claim. Similarly, the optimal local control in the
boundary between the non-action region (white) and the dividend payment region
for the first company (dark grey region), would be (assuming some regularity
on the boundary) that the first company pay dividends with some rate and the
control surplus would move leftward until the bivariate surplus reaches the
point $(1.61,1.06)$.%

\[%
\begin{array}
[c]{ccc}%
{\parbox[b]{2.7588in}{\begin{center}
\includegraphics[
height=2.8072in,
width=2.7588in
]%
{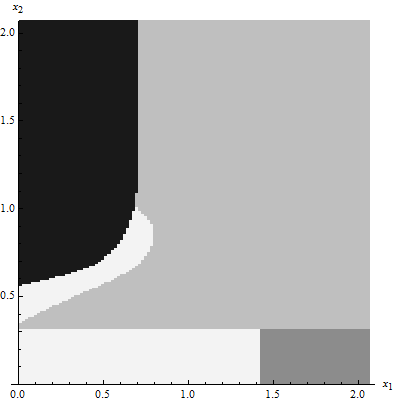}%
\\
Figure 1
\end{center}}}
&  &
{\parbox[b]{2.7501in}{\begin{center}
\includegraphics[
height=2.7985in,
width=2.7501in
]%
{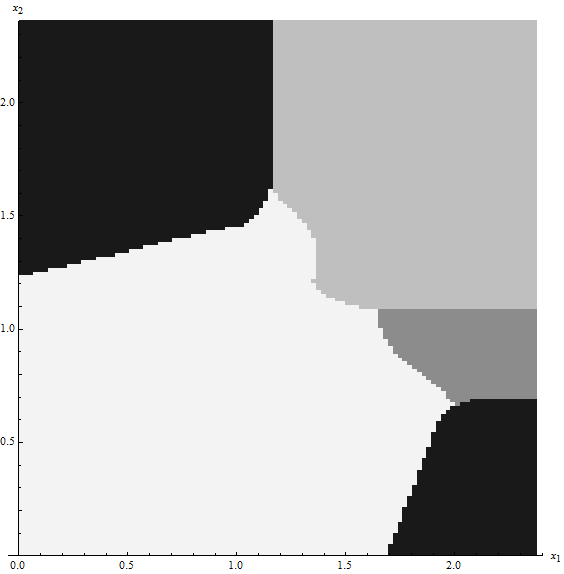}%
\\
Figure 2
\end{center}}}
\end{array}
\]

\section{Appendix}

This section contains the proofs of all the lemmas.

\textit{Proof of Lemma \ref{Dynkins}}. Let us extend the function $g$ to
$\mathbf{R}^{n}$ as $g(\mathbf{x})=0$ for $\mathbf{x}\notin\mathbf{R}_{+}^{n}
$ and the function $\upsilon\ $to $\mathbf{R}^{n}\times\mathbf{R}_{+}^{n}$ as
$\upsilon(\mathbf{x},\mathbf{\alpha})=0$ for $\left(  \mathbf{R}^{n}%
\times\mathbf{R}_{+}^{n}\right)  \diagup B$, where $B$ is defined in
(\ref{Definicion B}). Using the expressions (\ref{UncontrolledSurplusOriginal}%
) and the change of variables formula for finite variation processes, and
calling $\mathbf{z}_{s}=\mathbf{X}_{s}^{\mathbf{L}}{}$ and $\mathbf{\breve{z}%
}_{s}=\mathbf{\check{X}}_{s}^{\mathbf{L}}$, we can write%
\begin{equation}%
\begin{array}
[c]{l}%
g(\mathbf{z}_{\tau})e^{-c\tau}-g(\mathbf{x})\\%
\begin{array}
[c]{ll}%
= & \int\nolimits_{0}^{\tau}\mathbf{p\cdot}\nabla g(\mathbf{z}_{s^{-}}%
)e^{-cs}ds-c\int\nolimits_{0}^{\tau}g(\mathbf{z}_{s^{-}})e^{-cs}ds\\
& -\int\nolimits_{0}^{\tau}e^{-cs}\left(  \nabla g(\mathbf{z}_{s^{-}%
})\mathbf{\cdot}d\mathbf{L}_{s}^{c}\right)  +\sum\limits_{\mathbf{L}_{s}%
\neq\mathbf{L}_{s^{-}},~s\leq\tau}\left(  g(\mathbf{z}_{s})-g(\mathbf{\breve
{z}}_{s})\right)  e^{-cs}\\
& +\sum\limits_{\mathbf{\breve{z}}_{s}\neq\mathbf{z}_{s^{-}},~s\leq\tau
}\left(  g(\mathbf{\breve{z}}_{s})-g(\mathbf{z}_{s^{-}})\right)  e^{-cs}.
\end{array}
\end{array}
\label{Paso 1}%
\end{equation}
Note that $\mathbf{z}_{s}\in\mathbf{R}_{+}^{n}$ for $s\leq\tau$ except in the
case that $\tau=\tau^{\mathbf{L}}$. Since $\mathbf{z}_{s}=$ $\mathbf{\breve
{z}}_{s}-\Delta\mathbf{L}_{s},$%
\begin{equation}%
\begin{array}
[c]{l}%
-\int\nolimits_{0}^{\tau}e^{-cs}\nabla g(\mathbf{z}_{s^{-}})\mathbf{\cdot
}d\mathbf{L}_{s}^{c}+\sum\limits_{\mathbf{L}_{s}\neq\mathbf{L}_{s^{-}}%
,s\leq\tau}\left(  g(\mathbf{z}_{s})-g(\mathbf{\breve{z}}_{s})\right)
e^{-cs}\\%
\begin{array}
[c]{cl}%
= & -\int\nolimits_{0}^{\tau}e^{-cs}\nabla g(\mathbf{z}_{s^{-}})\mathbf{\cdot
}d\mathbf{L}_{s}^{c}-\sum\limits_{\mathbf{L}_{s}\neq\mathbf{L}_{s^{-}}%
,s\leq\tau}e^{-cs}\left(  \int\nolimits_{0}^{1}\left(  \nabla g\left(
\mathbf{\breve{z}}_{s}-\gamma\Delta\mathbf{L}_{s}\right)  \mathbf{\cdot}%
\Delta\mathbf{L}_{s}\right)  d\gamma\right) \\
= & -\int_{0^{-}}^{\tau}e^{-cs}\mathbf{a}\cdot d\mathbf{L}_{s}+\int
\nolimits_{0}^{\tau}e^{-cs}\left(  \mathbf{a}-\nabla g(\mathbf{z}_{s^{-}%
})\right)  \mathbf{\cdot}d\mathbf{L}_{s}^{c}\\
& +\sum\limits_{\mathbf{L}_{s}\neq\mathbf{L}_{s^{-}},s\leq\tau}e^{-cs}%
\int\nolimits_{0}^{1}\left(  \mathbf{a}-\nabla g\left(  \mathbf{\breve{z}}%
_{s}-\gamma\Delta\mathbf{L}_{s}\right)  \right)  \mathbf{\cdot}\Delta
\mathbf{L}_{s}d\gamma.
\end{array}
\end{array}
\label{Paso 2}%
\end{equation}
Since
\begin{equation}
M_{1}(t)=\sum\limits_{\mathbf{\breve{z}}\left(  s^{-}\right)  \neq
\mathbf{z}_{s^{-}},s\leq t}\left(  g(\mathbf{\breve{z}}_{s})-g(\mathbf{z}%
_{s^{-}})\right)  e^{-cs}-\lambda\int\limits_{0}^{t}e^{-cs}\int
\limits_{\mathbf{R}_{+}^{n}}\left(  g(\mathbf{z}_{s^{-}}-\mathbf{\alpha
})-g(\mathbf{z}_{s^{-}})\right)  dF(\mathbf{\alpha})ds\label{M1}%
\end{equation}
and%
\begin{equation}
M_{2}(t)=\sum\limits_{\mathbf{\breve{z}}\left(  s^{-}\right)  \neq
\mathbf{z}_{s^{-}},s\leq t}-\upsilon(\mathbf{\breve{z}}_{s^{-}},\mathbf{z}%
(s^{-})-\mathbf{\breve{z}}_{s})e^{-cs}+\lambda\int\limits_{0}^{t}e^{-cs}%
\int\limits_{\mathbf{R}_{+}^{n}}\upsilon(\mathbf{z}_{s^{-}},\mathbf{\alpha
})dF(\mathbf{\alpha})ds\label{M2}%
\end{equation}
are martingales with zero expectation, we have from (\ref{Paso 1}) and
(\ref{Paso 2})%
\[%
\begin{array}
[c]{l}%
(g(\mathbf{z}_{\tau})I_{\{\tau<\tau^{\mathbf{L}}\}}-\upsilon(\mathbf{z}%
_{\tau^{-}},\mathbf{z}_{\tau^{-}}-\mathbf{z}_{\tau})I_{\{\tau=\tau
^{\mathbf{L}}\}})e^{-c\tau}-g(\mathbf{x})\\%
\begin{array}
[c]{ll}%
= & (g(\mathbf{z}_{\tau})-\upsilon(\mathbf{z}_{\tau^{-}},\mathbf{z}_{\tau^{-}%
}-\mathbf{z}(\tau)))e^{-c\tau}-g(\mathbf{x})\\
= & \int\nolimits_{0}^{\tau}\mathcal{L}(g)(\mathbf{z}_{s^{-}})e^{-cs}%
ds-\int_{0^{-}}^{\tau}e^{-cs}\mathbf{a}\cdot d\mathbf{L}_{s}\\
& +\int\nolimits_{0}^{\tau}e^{-cs}\left(  \mathbf{a}-\nabla g(\mathbf{z}%
_{s^{-}})\right)  \mathbf{\cdot}d\mathbf{L}_{s}^{c}\\
& +\sum\limits_{\mathbf{L}_{s}\neq\mathbf{L}_{s^{-}},s\leq\tau}e^{-cs}%
\int\nolimits_{0}^{1}\left(  \mathbf{a}-\nabla g\left(  \mathbf{\breve{z}}%
_{s}-\gamma\Delta\mathbf{L}_{s}\right)  \mathbf{\cdot}\Delta\mathbf{L}%
_{s}\right)  d\gamma+M(\tau);
\end{array}
\end{array}
\]
where $M(t)=M_{1}(t)+M_{2}(t)$. $\blacksquare$

In order to prove Lemma \ref{SupersolucionMayor-ValueFunction}, we will use a
technical lemma in which we construct a sequence of smooth functions that
approximate a (possible non-smooth) viscosity supersolution. This is done in
order to apply Lemma \ref{Dynkins} to an approximate smooth function instead
of the viscosity supersolution; we have to do that because the amount of time
the controlled process spends at non-differentiable points of the viscosity
supersolution could have positive Lebesgue measure. We omit the proof of this
lemma because it is similar to the one-dimensional version given in Lemma 4.1
of \cite{AM Libro}; the result is obtained by standard convolution arguments
using that the function $\mathcal{R}$ is continuous.

\begin{lemma}
\label{A.1} Fix $\mathbf{x}^{0}\ $in the interior of $\mathbf{R}_{+}^{n}$ and
let $\overline{u}$\ be a supersolution of (\ref{HJB}) satisfying the growth
condition (\ref{gc}).\ We can find a sequence of functions $\overline{u}%
_{m}:\mathbf{R}_{+}^{n}\rightarrow\mathbf{R}$\ such that:

(a) $\overline{u}_{m}$\ is continuously differentiable and $\overline{u}%
_{m}\geq\overline{u}\geq f.$

(b) $\overline{u}_{m}\ $satisfies the growth condition (\ref{gc}).

(c)$\ \mathbf{p\cdot}\nabla\overline{u}_{m}$\ $\leq\left(  c+\lambda\right)
\overline{u}_{m}+\lambda\left\vert \overline{u}(\mathbf{0})\right\vert
+\lambda\mathbb{E}\left(  \left\vert \upsilon(\mathbf{0},\mathbf{U}%
_{1})\right\vert \right)  $ in $\mathbf{R}_{+}^{n}$ and $\mathbf{a}%
-\nabla\overline{u}_{m}\leq\mathbf{0}$.

(d) $\overline{u}_{m}$\ $\searrow$ $\overline{u}$\ uniformly on compact sets
in $\mathbf{R}_{+}^{n}$ and $\nabla\overline{u}_{m}$\ converges to
$\nabla\overline{u}$\ a.e. in $\mathbf{R}_{+}^{n}$.

(e) There exists a sequence $c_{m}$ with $\lim\limits_{m\rightarrow\infty
}c_{m}=0$ such that
\[
\sup\nolimits_{\mathbf{x}\in\lbrack\mathbf{0},\mathbf{x}^{0}]}\mathcal{L}%
(\overline{u}_{m})\left(  \mathbf{x}\right)  \leq c_{m}.
\]

\textit{Proof of Lemma \ref{SupersolucionMayor-ValueFunction}}. Consider the
processes $\mathbf{z}_{s}=\mathbf{X}_{s}^{\mathbf{L}}{}$ defined in
(\ref{XL}), let us call $\tau=\tau^{\mathbf{L}}$ and take $\widetilde{\tau
}=\overline{\tau}\wedge\tau$. Let us consider the functions $\overline{u}_{m}$
defined in Lemma \ref{A.1} in $\mathbf{R}_{+}^{n}$ . Using Lemma \ref{Dynkins}
for $\widetilde{\tau}\wedge t$, we get from Lemma \ref{A.1} (a) and (c) that%
\begin{equation}%
\begin{array}
[c]{l}%
\overline{u}_{m}(\mathbf{z}_{t})e^{-ct}I_{\{t<\widetilde{\tau}\}}%
+e^{-c\overline{\tau}}f(\mathbf{z}_{\overline{\tau}})I_{\{t\wedge
\widetilde{\tau}=\overline{\tau},\overline{\tau}<\tau\}}-e^{-c\overline{\tau}%
}\upsilon\left(  \mathbf{z}_{\tau^{-}},\mathbf{z}_{\tau^{-}}-\mathbf{z}_{\tau
}\right)  I_{\{t\wedge\widetilde{\tau}=\tau\}}-\overline{u}_{m}(\mathbf{x})\\%
\begin{array}
[c]{ll}%
\leq & \overline{u}_{m}(\mathbf{z}_{t})e^{-ct}I_{\{t<\widetilde{\tau}%
\}}+e^{-c\overline{\tau}}\overline{u}_{m}(\mathbf{z}_{\overline{\tau}%
})I_{\{t\wedge\widetilde{\tau}=\overline{\tau},\overline{\tau}<\tau
\}}-e^{-c\overline{\tau}}\upsilon\left(  \mathbf{z}_{\tau^{-}},\mathbf{z}%
_{\tau^{-}}-\mathbf{z}_{\tau}\right)  I_{\{t\wedge\widetilde{\tau}=\tau
\}}-\overline{u}_{m}(\mathbf{x})\\
\leq & \int\nolimits_{0}^{t\wedge\widetilde{\tau}}\mathcal{L}(\overline{u}%
_{m})(\mathbf{z}_{s^{-}})e^{-cs}ds-\int_{0^{-}}^{t\wedge\widetilde{\tau}%
}e^{-cs}\mathbf{a}\cdot d\mathbf{L}_{s}+M(t\wedge\widetilde{\tau}),
\end{array}
\end{array}
\label{ItoUnMenorSuper}%
\end{equation}
where $M(t)\ $is a zero-expectation martingale. Since $\mathbf{L}_{s}$ is
non-decreasing we get, using the monotone convergence theorem, that%
\[%
\begin{array}
[c]{l}%
\lim\limits_{t\rightarrow\infty}\mathbb{E}_{\mathbf{x}}\left(  \int_{0^{-}%
}^{t\wedge\widetilde{\tau}}e^{-cs}\mathbf{a}\cdot d\mathbf{L}_{s}%
+e^{-c\overline{\tau}}f(\mathbf{z}_{\overline{\tau}})I_{\{t\wedge
\widetilde{\tau}=\overline{\tau},\overline{\tau}<\tau\}}-e^{-c\overline{\tau}%
}\upsilon\left(  \mathbf{z}_{\tau^{-}},\mathbf{z}_{\tau^{-}}-\mathbf{z}_{\tau
}\right)  I_{\{t\wedge\widetilde{\tau}=\tau\}}\right) \\
=V_{\pi}(\mathbf{x}).
\end{array}
\]
From Lemma \ref{A.1}-(c), we have%
\begin{equation}
-\left(  c+\lambda\right)  \overline{u}_{m}(\mathbf{x})+\overline{u}%
_{m}(0)\lambda F(\mathbf{x})-\lambda\mathbb{E}\left(  \left\vert
\upsilon(\mathbf{0},\mathbf{U}_{1})\right\vert \right)  \leq\mathcal{L}%
(\overline{u}_{m})(\mathbf{x})\leq\lambda\overline{u}_{m}(\mathbf{x}%
)+\lambda\left\vert \overline{u}(\mathbf{0})\right\vert +\lambda
\mathbb{E}\left(  \left\vert \upsilon(\mathbf{0},\mathbf{U}_{1})\right\vert
\right)  -\mathcal{R}(\mathbf{x}).\label{PrimeraCotaL}%
\end{equation}
By Lemma \ref{A.1}-(b), (c) and the inequality $\mathbf{z}_{s}\leq
\mathbf{x}+\mathbf{p}s,$ there exists $d_{0}$ large enough such that
\begin{equation}
\overline{u}_{m}(\mathbf{z}_{s})\leq\overline{u}_{m}(\mathbf{x}+\mathbf{p}%
s)\leq d_{0}e^{\frac{c}{2n}\sum_{i=1}^{n}\frac{x_{i}+p_{i}s}{p_{i}}}%
=d_{0}h_{0}(\mathbf{x})e^{\frac{c}{2}s}\label{Acotacionu}%
\end{equation}
and%
\begin{equation}
-\upsilon(\mathbf{z}_{s^{-}},\mathbf{\alpha})\leq S(\mathbf{z}_{s^{-}})\leq
d_{0}h_{0}(\mathbf{x})e^{\frac{c}{2}s}\text{ for }\left(  \mathbf{z}_{s^{-}%
}-\mathbf{\alpha}\right)  \notin\mathbf{R}_{+}^{n},\label{Acotacion_Rho}%
\end{equation}
where $h_{0}$ and $S$ are defined in (\ref{ho}) and Proposition
\ref{Crecimiento de V} respectively. Therefore, from (\ref{PrimeraCotaL}), we
obtain that there exists $d_{1}$ large enough such that,
\begin{equation}
e^{-cs}\left\vert \mathcal{L}(\overline{u}_{m})\left(  \mathbf{z}_{s^{-}%
}\right)  \right\vert \leq d_{1}e^{-\frac{c}{2}s}.\label{Acotacion_L}%
\end{equation}
And using the bounded convergence theorem,%
\begin{equation}
\lim\limits_{t\rightarrow\infty}\mathbb{E}_{\mathbf{x}}\left(  \int
\nolimits_{0}^{t\wedge\widetilde{\tau}}\mathcal{L}(\overline{u}_{m}%
)(\mathbf{z}_{s^{-}})e^{-cs}ds\right)  =\mathbb{E}_{\mathbf{x}}\left(
\int\nolimits_{0}^{\widetilde{\tau}}\mathcal{L}(\overline{u}_{m}%
)(\mathbf{z}_{s^{-}})e^{-cs}ds\right)  .\label{monotone2}%
\end{equation}
From (\ref{ItoUnMenorSuper}) and (\ref{monotone2}), we get
\begin{equation}
\lim\limits_{t\rightarrow\infty}\mathbb{E}_{\mathbf{x}}\left(  \overline
{u}_{m}(\mathbf{z}_{t})e^{-ct}I_{\{t<\widetilde{\tau}\}}\right)  -\overline
{u}_{m}(\mathbf{x})\leq\mathbb{E}_{\mathbf{x}}\left(  \int\nolimits_{0}%
^{\widetilde{\tau}}\mathcal{L}(\overline{u}_{m})(\mathbf{z}_{s^{-}}%
)e^{-cs}ds\right)  -V_{\pi}(\mathbf{x}).\label{limite0}%
\end{equation}
By (\ref{Acotacionu}),
\begin{equation}
\lim\limits_{t\rightarrow\infty}\mathbb{E}_{\mathbf{x}}\left(  \overline
{u}_{m}(\mathbf{z}_{t})e^{-ct}I_{\{t<\widetilde{\tau}\}}\right)
=0.\label{limite1}%
\end{equation}
Let us prove now that
\begin{equation}
\limsup\limits_{m\rightarrow\infty}\mathbb{E}_{\mathbf{x}}\left(
\int\nolimits_{0}^{\widetilde{\tau}}\mathcal{L}(\overline{u}_{m}%
)(\mathbf{z}_{s^{-}})e^{-cs}ds\right)  \leq0.\label{limite2}%
\end{equation}
Given any $\varepsilon>0$, from (\ref{Acotacion_L}), we can find $T$ large
enough such that
\begin{equation}
\mathbb{E}_{\mathbf{x}}\left(  \int\nolimits_{T\wedge\widetilde{\tau}%
}^{\widetilde{\tau}}\left\vert \mathcal{L}(\overline{u}_{m})(\mathbf{z}%
_{s^{-}})\right\vert e^{-cs}ds\right)  \leq\frac{2d_{1}}{c}(e^{-\frac{c}{2}%
T})<\frac{\varepsilon}{2}.\label{arreglado5}%
\end{equation}
For $s\leq T$, we get $\mathbf{z}_{s^{-}}\in\lbrack\mathbf{0},$ $\mathbf{x}%
+\mathbf{p}T$ $]$ , then from Lemma \ref{A.1}-(e) we can find $m_{0}$ large
enough such that for any $m\geq m_{0}$%
\[
\int\nolimits_{0}^{T}\mathcal{L}(\overline{u}_{m})(\mathbf{z}_{s^{-}}%
)e^{-cs}ds\leq c_{m}\int\nolimits_{0}^{T}e^{-cs}ds\leq\frac{c_{m}}{c}\leq
\frac{\varepsilon}{2}%
\]
and so we have (\ref{limite2}). Thus, from (\ref{limite0}) and using
(\ref{limite1}) and (\ref{limite2}), we obtain
\end{lemma}%

\begin{equation}
\overline{u}(\mathbf{x})=\lim\nolimits_{m\rightarrow\infty}\overline{u}%
_{m}(\mathbf{x})\geq V_{\pi}(\mathbf{x})\text{. }\blacksquare
\label{supersolucionMayorqueVL}%
\end{equation}

\textit{Proof of Lemma \ref{Lema de tiempo infinito}. }Suppose that$\ \tilde
{k}=\infty$, calling
\[
k_{l}:=\mathbf{m\cdot1}+(l-1)n+1,
\]
there are at least $i_{l}\geq l$ control actions $\mathbf{E}_{0}$ in $\left(
s_{1},s_{2},...,s_{k_{l}}\right)  $. Let us consider the non-decreasing
sequence $(j_{l})_{l}$ defined as
\[
j_{l}:=\max\{j:\tau_{j}\leq t_{k_{l}}\},
\]
we have that $t_{k_{l}}\geq\tau_{j_{l}}+(i_{l}-j_{l})\delta$. If
$\lim_{l\rightarrow\infty}i_{l}-j_{l}=\infty$, then
\[
\lim\nolimits_{l\rightarrow\infty}t_{k_{l}}\geq\lim\nolimits_{l\rightarrow
\infty}\tau_{j_{l}}+(i_{l}-j_{l})\delta\geq\lim\nolimits_{l\rightarrow\infty
}(i_{l}-j_{l})\delta=\infty;
\]
if not, $\lim_{l\rightarrow\infty}j_{l}=\infty$ and so%

\[
\lim\nolimits_{l\rightarrow\infty}t_{k_{l}}\geq\lim\nolimits_{l\rightarrow
\infty}\tau_{j_{l}}+(i_{l}-j_{l})\delta\geq\lim\nolimits_{l\rightarrow\infty
}\tau_{j_{l}}%
\]
and since $\lim_{l\rightarrow\infty}\tau_{j_{l}}=$ $\lim_{i\rightarrow\infty
}\tau_{i}=$ $\infty$ a.s., we have the result. $\blacksquare$

\textit{Proof of Lemma \ref{Ts crecientes}. }It is straightforward
that\textit{\ }$T_{0}$, $T_{i},$ $T_{s}$ and $T$ are non-decreasing and that
\[
\sup\nolimits_{\mathbf{m}\in\mathbf{N}_{0}^{n}}\left\vert T(w_{1}%
)(\mathbf{m})-T(w_{2})(\mathbf{m})\right\vert \leq\sup\nolimits_{\mathbf{m}%
\in\mathbf{N}_{0}^{n}}\left\vert w_{1}(\mathbf{m})-w_{2}(\mathbf{m}%
)\right\vert .
\]

Also, given a function $w:\mathbf{N}_{0}^{n}\rightarrow\mathbf{R}$ it is
immediate to see that $T_{i}(w)$ and $T_{s}(w)$ can be written as a linear
combination of the values of $w(\mathbf{m})$ plus a constant. Let us prove now
that
\[
T_{0}(w)(\mathbf{m})=e^{-(c+\lambda)\delta}w(\mathbf{m}+\mathbf{1}%
)+\sum\nolimits_{0\leq\mathbf{k}\leq\mathbf{m}}a_{1}(\mathbf{k},\mathbf{m}%
)w(\mathbf{k})+a_{2}(\mathbf{m})\text{,}%
\]

\begin{lemma}
where%
\[%
\begin{array}
[c]{lll}%
a_{1}(\mathbf{k},\mathbf{m}) & = & I_{\{\mathbf{k}\leq\mathbf{m}-\mathbf{1}%
\}}\int\limits_{0}^{\delta}\lambda e^{-(c+\lambda)t}(F(g^{\delta}\left(
\mathbf{m}-\mathbf{k}\right)  +t\mathbf{p})-F(g^{\delta}\left(  \mathbf{m}%
-\mathbf{k}-\mathbf{1}\right)  +t\mathbf{p}))dt\text{ \ }\\
&  & +I_{\{\mathbf{k}\leq\mathbf{m},\mathbf{k}\nleqslant\mathbf{m}%
-\mathbf{1}\}}\int\limits_{0}^{\delta}\lambda e^{-(c+\lambda)t}(F(g^{\delta
}\left(  \mathbf{m}-\mathbf{k}\right)  +t\mathbf{p})-F(\mathbf{0}\vee\left(
g^{\delta}\left(  \mathbf{m}-\mathbf{k}\right)  +t\mathbf{p}\right)  ))dt~
\end{array}
\]

and%
\[%
\begin{array}
[c]{lll}%
a_{2}(\mathbf{m}) & = & \sum\limits_{0\leq\mathbf{k}<\mathbf{m}-\mathbf{1}%
}\int\limits_{0}^{\delta}(\lambda e^{-(c+\lambda)t}%
{\textstyle\int\limits_{g^{\delta}\left(  \mathbf{m}-\mathbf{k}-\mathbf{1}%
\right)  +t\mathbf{p}}^{g^{\delta}\left(  \mathbf{m}-\mathbf{k}\right)
+t\mathbf{p}}}
\mathbf{a}\cdot(g^{\delta}\left(  \mathbf{m}-\mathbf{k}\right)  +t\mathbf{p}%
-\mathbf{\alpha})dF(\mathbf{\alpha}))dt\\
&  & +\sum\limits_{\mathbf{k}\leq\mathbf{m},\mathbf{k}\nleqslant
\mathbf{m}-\mathbf{1}}\int\limits_{0}^{\delta}(\lambda e^{-(c+\lambda)t}%
{\textstyle\int\limits_{\mathbf{0}\vee\left(  g^{\delta}\left(  \mathbf{m}%
-\mathbf{k}\right)  +t\mathbf{p}\right)  }^{g^{\delta}\left(  \mathbf{m}%
-\mathbf{k}\right)  +t\mathbf{p}}}
\mathbf{a}\cdot(g^{\delta}\left(  \mathbf{m}-\mathbf{k}\right)  +t\mathbf{p}%
-\mathbf{\alpha})dF(\mathbf{\alpha}))dt\\
&  & -\int\limits_{0}^{\delta}e^{-(c+\lambda)t}\mathcal{R}(g^{\delta
}(\mathbf{m})+t\mathbf{p})dt.
\end{array}
\]

\end{lemma}

\textit{\ }Given $\mathbf{m}\in\mathbf{N}_{0}^{n}$, $\mathbf{\alpha}%
\in\mathbf{R}_{+}^{n}$ and $0<t\leq\delta$ such that $\mathbf{0}\leq
g^{\delta}(\mathbf{m})+t\mathbf{p}-\mathbf{\alpha}$, let us define
\[
\mathbf{k}:=\rho^{\delta}(g^{\delta}(\mathbf{m})+t\mathbf{p}-\mathbf{\alpha}),
\]
and so $\mathbf{k}\leq\mathbf{m}.$

If $\mathbf{k}\leq\mathbf{m}-\mathbf{1}$,
\[
g^{\delta}(\mathbf{k})\leq g^{\delta}(\mathbf{m})+t\mathbf{p}-\mathbf{\alpha
}<g^{\delta}\left(  \mathbf{k}+\mathbf{1}\right)  \leq g^{\delta}(\mathbf{m})
\]
that implies
\[
\mathbf{0}<g^{\delta}\left(  \mathbf{m}-\mathbf{k}-\mathbf{1}\right)
+t\mathbf{p}<\mathbf{\alpha}\leq g^{\delta}\left(  \mathbf{m}-\mathbf{k}%
\right)  +t\mathbf{p}.
\]

If $\mathbf{k}\leq\mathbf{m}$ with $\mathbf{k}\nleqslant\mathbf{m}-\mathbf{1}%
$,%
\[
g^{\delta}(\mathbf{k})\leq g^{\delta}(\mathbf{m})+t\mathbf{p}-\mathbf{\alpha
}<g^{\delta}\left(  \mathbf{k}+\mathbf{1}\right)  \wedge\left(  g^{\delta
}(\mathbf{m})+t\mathbf{p}\right)
\]
and so%
\[
\left(  g^{\delta}\left(  \mathbf{m}-\mathbf{k}-\mathbf{1}\right)
+t\mathbf{p}\right)  \vee\mathbf{0}<\mathbf{\alpha}\leq g^{\delta}\left(
\mathbf{m}-\mathbf{k}\right)  +t\mathbf{p}.
\]
Then, we can write
\[%
\begin{array}
[c]{l}%
\mathcal{I}^{\delta}(w)(\mathbf{m})\\%
\begin{array}
[c]{ll}%
= & \sum_{0\leq\mathbf{k}\leq\mathbf{m}-\mathbf{1}}w(\mathbf{k})\int
_{0}^{\delta}\lambda e^{-(c+\lambda)t}(%
{\textstyle\int\nolimits_{g^{\delta}\left(  \mathbf{m}-\mathbf{k}%
-\mathbf{1}\right)  +t\mathbf{p}}^{g^{\delta}\left(  \mathbf{m}-\mathbf{k}%
\right)  +t\mathbf{p}}}
dF(\mathbf{\alpha}))dt\\
& +\sum_{0\leq\mathbf{k}\leq\mathbf{m}-\mathbf{1}}\int_{0}^{\delta}\lambda
e^{-(c+\lambda)t}(%
{\textstyle\int\nolimits_{g^{\delta}\left(  \mathbf{m}-\mathbf{k}%
-\mathbf{1}\right)  +t\mathbf{p}}^{g^{\delta}\left(  \mathbf{m}-\mathbf{k}%
\right)  +t\mathbf{p}}}
\mathbf{a}\cdot\left(  g^{\delta}(\mathbf{m}-\mathbf{k})+t\mathbf{p}%
-\mathbf{\alpha}\right)  dF(\mathbf{\alpha}))dt\\
& +\sum_{\mathbf{k}\leq\mathbf{m},\mathbf{k}\nleqslant\mathbf{m}-\mathbf{1}%
}w(\mathbf{k})\int_{0}^{\delta}\lambda e^{-(c+\lambda)t}(%
{\textstyle\int\nolimits_{\left(  g^{\delta}\left(  \mathbf{m}-\mathbf{k}%
-\mathbf{1}\right)  +t\mathbf{p}\right)  \vee\mathbf{0}}^{g^{\delta}\left(
\mathbf{m}-\mathbf{k}\right)  +t\mathbf{p}}}
dF(\mathbf{\alpha}))dt\\
& +\sum_{\mathbf{k}\leq\mathbf{m},\mathbf{k}\nleqslant\mathbf{m}-\mathbf{1}%
}\int_{0}^{\delta}\lambda e^{-(c+\lambda)t}(%
{\textstyle\int\nolimits_{\left(  g^{\delta}\left(  \mathbf{m}-\mathbf{k}%
-\mathbf{1}\right)  +t\mathbf{p}\right)  \vee\mathbf{0}}^{g^{\delta}\left(
\mathbf{m}-\mathbf{k}\right)  +t\mathbf{p}}}
\mathbf{a}\cdot\left(  g^{\delta}(\mathbf{m}-\mathbf{k})+t\mathbf{p}%
-\mathbf{\alpha}\right)  dF(\mathbf{\alpha}))dt.
\end{array}
\end{array}
\]
Therefore, from (\ref{Definicion TE}), we have the result. $\blacksquare$

\textit{Proof of Lemma \ref{Menor Supersolucion Discreta}}. The proof of this
lemma is a discrete version of the one of Lemma
\ref{SupersolucionMayor-ValueFunction}.\textit{\ }Assume that $\pi
=(\mathbf{L},\overline{\tau})\in\Pi_{g^{\delta}(\mathbf{m})}^{\delta}$. For
any $\omega=(\tau_{i},\mathbf{U}_{i})_{i\geq1}$, consider the sequence
$\mathbf{s}=(s_{k})_{k=1,...,\tilde{k}}$ with $s_{k}\in\mathcal{E}$
corresponding to $\pi$ and $\mathbf{m}^{k}$, $\mathbf{y}^{k}$ and times
$t_{k}$ and $\Delta_{k}$ as defined in Section \ref{Discrete Scheme}. Let
$\left(  \kappa_{l}\right)  _{l\geq1}$ be the indices of the sequence
$\mathbf{s}=(s_{k})_{k=1,...,\tilde{k}}$ where $s_{k}$ is either
$\mathbf{E}_{s}$ or $\mathbf{E}_{0}$\textbf{. }If the sequence stops at
$\tilde{k}=\kappa_{l_{0}}<\infty$, we define
\[
\kappa_{l}=\kappa_{l_{0}}\text{ for }l\geq l_{0},\text{ }t_{\kappa_{l_{0}+j}%
}=t_{\kappa_{l_{0}}}+\Delta_{\kappa_{l_{0}}}\text{ for }j\geq1;
\]
and if $\tilde{k}=\infty$ we put $l_{0}=\infty$. Consider the case in which
the process goes to ruin at $\kappa_{l}$, that is $\mathbf{y}^{\kappa_{l}%
}\notin\mathbf{R}_{+}^{n}$; then the surplus prior to the ruin is
$\mathbf{y}^{\kappa_{l}}+\mathbf{U}$ and the penalty paid at ruin is
$\upsilon(\mathbf{y}^{\kappa_{l}}+\mathbf{U},\mathbf{U})$, where $\mathbf{U} $
is the last jump of the uncontrolled process. So we define, for $l\geq1$,
\[
H(l)=w(\mathbf{m}^{1+\kappa_{l}})I_{\{s_{\kappa_{l}}=\mathbf{E}_{0}%
\}}I_{\{\mathbf{y}^{\kappa_{l}}\in\mathbf{R}_{+}^{n}\}}-\upsilon
(\mathbf{y}^{\kappa_{l}}+\mathbf{U},\mathbf{U})I_{\{s_{\kappa_{l}}%
=\mathbf{E}_{0}\}}I_{\{\mathbf{y}^{\kappa_{l}}\notin\mathbf{R}_{+}^{n}%
\}}+f(g^{\delta}\left(  \mathbf{m}^{\kappa_{l}}\right)  )I_{\{s_{\kappa_{l}%
}=\mathbf{E}_{s}\}}\text{.}%
\]
If we put $H(0)=w(\mathbf{m})$, $\kappa_{0}=0$ and $t_{0}=0$, we have using
$\left(  T_{i}(w)-w\right)  _{i=1,...,n}\leq0,$%
\begin{equation}%
\begin{array}
[c]{lll}%
e^{-ct_{\kappa_{l+1}}}H(l)-w(\mathbf{m}) & = & \sum_{j=1}^{l}(e^{-ct_{\kappa
_{j+1}}}H(j)-e^{-ct_{\kappa_{j}}}H(j-1))\\
& = & \sum_{j=1}^{l}I_{\{\kappa_{j+1}\neq\kappa_{j}\}}(e^{-ct_{\kappa_{j+1}}%
}H(j)-e^{-ct_{\kappa_{j}}}H(j-1))\\
& = & \sum_{j=1}^{l}I_{\{\kappa_{j+1}\neq\kappa_{j}\}}(e^{-ct_{1+\kappa_{j-1}%
}}(\sum_{k=1+\kappa_{j-1}}^{\kappa_{j}-1}\left(  w(\mathbf{m}^{k+1}%
)-w(\mathbf{m}^{k})\right)  ))\\
&  & +\sum_{j=1}^{l}I_{\{\kappa_{j+1}\neq\kappa_{j}\}}(e^{-ct_{\kappa_{j+1}}%
}H(j)-e^{-ct_{\kappa_{j}}}w(\mathbf{m}^{\kappa_{j}}))\\
& \leq & \sum_{j=1}^{l}I_{\{\kappa_{j+1}\neq\kappa_{j}\}}(\sum_{k=1+\kappa
_{j-1}}^{\kappa_{j}-1}e^{-ct_{1+\kappa_{j-1}}}(\sum_{i=1}^{n}\left(
-a_{i}p_{i}\delta\right)  I_{\{s_{k}=\mathbf{E}_{i}\}}))\\
&  & +\sum_{j=1}^{l}I_{\{\kappa_{j+1}\neq\kappa_{j}\}}(e^{-ct_{\kappa_{j+1}}%
}H(j)-e^{-ct_{\kappa_{j}}}w(\mathbf{m}^{\kappa_{j}}));
\end{array}
\label{Desigualdad 1}%
\end{equation}
and since $T_{0}(w)-w\leq0$ and $T_{s}(w)-w\leq0,$ if $\kappa_{j+1}\neq
\kappa_{j}$,%
\begin{equation}%
\begin{array}
[c]{l}%
\mathbb{E}\left(  \left.  e^{-ct_{\kappa_{j+1}}}H(j)-e^{-ct_{\kappa_{j}}%
}w(\mathbf{m}^{\kappa_{j}})\right\vert \mathcal{F}_{t_{\kappa_{j}}}\right) \\%
\begin{array}
[c]{ll}%
= & \mathbb{E}\left(  \left.  (e^{-ct_{\kappa_{j+1}}}H(j)-e^{-ct_{\kappa_{j}}%
}w(\mathbf{m}^{\kappa_{j}}))I_{\{s_{\kappa_{j}}=\mathbf{E}_{0}\}}\right\vert
\mathcal{F}_{t_{\kappa_{j}}}\right)  +I_{\{s_{\kappa_{j}}=\mathbf{E}_{s}%
\}}e^{-ct_{\kappa_{j}}}\left(  f(g^{\delta}(\mathbf{m}^{\kappa_{j}%
}))-w(\mathbf{m}^{\kappa_{j}})\right) \\
\leq & \mathbb{E}\left(  \left.  e^{-ct_{\kappa_{j+1}}}I_{\{s_{\kappa_{j}%
}=\mathbf{E}_{0}\}}(w(\mathbf{m}^{1+\kappa_{j}})I_{\{\mathbf{y}^{\kappa_{j}%
}\in\mathbf{R}_{+}^{n}\}}-\upsilon(\mathbf{y}^{\kappa_{j}}+\mathbf{U}%
,\mathbf{U})I_{\{\mathbf{y}_{\kappa_{j}}\notin\mathbf{R}_{+}^{n}%
\}})\right\vert \mathcal{F}_{t_{\kappa_{j}}}\right) \\
& -e^{-ct_{\kappa_{j}}}w(\mathbf{m}^{\kappa_{j}})I_{\{s_{\kappa_{j}%
}=\mathbf{E}_{0}\}}\\
= & e^{-ct_{\kappa_{j}}}I_{\{s_{\kappa_{j}}=\mathbf{E}_{0}\}}\left(
T_{0}(w)\left(  \mathbf{m}^{\kappa_{j}}\right)  -w(\mathbf{m}^{\kappa_{j}%
})\right) \\
& -e^{-ct_{\kappa_{j}}}I_{\{s_{\kappa_{j}}=\mathbf{E}_{0}\}}\int
\limits_{0}^{\delta}%
{\textstyle\int\limits_{\mathbf{\alpha}\in\lbrack\mathbf{0},\mathbf{z}%
_{j}(t)]}}
\lambda e^{-(c+\lambda)t}\mathbf{a}\cdot\left(  \mathbf{z}_{j}%
(t)-\mathbf{\alpha}-\left\langle \mathbf{z}_{j}(t)-\mathbf{\alpha
}\right\rangle ^{\delta}\right)  dF(\mathbf{\alpha})dt\\
\leq & -e^{-ct_{\kappa_{j}}}I_{\{s_{\kappa_{j}}=\mathbf{E}_{0}\}}%
\int\limits_{0}^{\delta}%
{\textstyle\int\limits_{\mathbf{\alpha}\in\lbrack\mathbf{0},\mathbf{z}%
_{j}(t)]}}
\lambda e^{-(c+\lambda)t}\mathbf{a}\cdot\left(  \mathbf{z}_{j}%
(t)-\mathbf{\alpha}-\left\langle \mathbf{z}_{j}(t)-\mathbf{\alpha
}\right\rangle ^{\delta}\right)  dF(\mathbf{\alpha})dt\text{,}%
\end{array}
\end{array}
\label{Desigualdad 2}%
\end{equation}
where $\mathbf{z}_{j}(t)=g^{\delta}(\mathbf{m}^{\kappa_{j}})+t\mathbf{p}$.
From (\ref{Desigualdad 1}) and (\ref{Desigualdad 2}), and calling the initial
surplus $\mathbf{x}=g^{\delta}(\mathbf{m})\in\mathcal{G}^{\delta}$ we have,
\[
\lim\sup_{l\rightarrow\infty}\mathbb{E}_{\mathbf{x}}\left(  e^{-ct_{\kappa
_{l+1}}}H(l)-w(\mathbf{m})\right)  \leq-\mathbb{E}_{\mathbf{x}}\left(
\int_{0^{-}}^{\overline{\tau}\wedge\tau_{L}}e^{-cs}\mathbf{a}\cdot
d\mathbf{L}_{s}\right)  .
\]
Then,%
\[
w(\mathbf{m})\geq V_{\pi}(g^{\delta}(\mathbf{m}))+\lim\sup_{l\rightarrow
\infty}\mathbb{E}_{\mathbf{x}}\left(  I_{\{l\leq l_{0}\}}e^{-ct_{1+\kappa_{l}%
}}w(\mathbf{m}^{1+\kappa_{l}})I_{\{\mathbf{y}^{\kappa_{l}}\in\mathbf{R}%
_{+}^{n}\}}\right)  .
\]
Since%
\[
g^{\delta}(\mathbf{m}^{1+\kappa_{l}})\leq g^{\delta}\left(  \mathbf{m}%
+\rho^{\delta}(t_{1+\kappa_{l}}\mathbf{p})\right)
\]
and $w$ satisfies the growth condition (\ref{Growth Condition discreta}),
there exists $d$ large enough such that
\[
\lim\sup_{l\rightarrow\infty}\left(  \mathbb{E}_{\mathbf{x}}I_{\{l\leq
l_{0}\}}e^{-ct_{1+\kappa_{l}}}w(\mathbf{m}^{1+\kappa_{l}})I_{\{\mathbf{y}%
^{\kappa_{l}}\in\mathbf{R}_{+}^{n}\}}\right)  \leq d\lim_{l\rightarrow\infty
}\mathbb{E}_{\mathbf{x}}\left(  I_{\{l\leq l_{0}\}}e^{-ct_{1+\kappa_{l}}%
}e^{c\delta\mathbf{m\cdot1}/\left(  2n\right)  }e^{\frac{c}{2}t_{1+\kappa_{l}%
}})\right)  =0;
\]
so we have the result. $\blacksquare$

\textit{Proof of Lemma \ref{v_delta es lipschitz y monotona en el reticulado}}.

(1) Take the $\mathcal{G}^{\delta}$-optimal strategy $\pi_{g^{\delta
}(\mathbf{m})}^{\delta}\in\Pi_{g^{\delta}(\mathbf{m})}^{\delta}$ and define
$\overline{\pi}_{g^{\delta}(\mathbf{m}+\mathbf{e}_{i})}\in\Pi_{g^{\delta
}(\mathbf{m}+\mathbf{e}_{i})}^{\delta}$ \ by applying first the control action
$\mathbf{E}_{i}$ and then the $\mathcal{G}^{\delta}$-optimal strategy
$\pi_{g^{\delta}(\mathbf{m})}^{\delta}$. The value function of this strategy
is given by
\[
a_{i}p_{i}\delta+v^{\delta}(\mathbf{m}),
\]
so we obtain the the first inequality of this proposition. Now, take the
$\mathcal{G}^{\delta}$-optimal strategy $\pi_{g^{\delta}\left(  \mathbf{m}%
+\mathbf{1}\right)  }^{\delta}\in\Pi_{g^{\delta}\left(  \mathbf{m}%
+\mathbf{1}\right)  }^{\delta}$ and define $\overline{\pi}_{g^{\delta
}(\mathbf{m})}\in\Pi_{g^{\delta}(\mathbf{m})}^{\delta}$ by applying first the
control action $\mathbf{E}_{0}$ and then the $\mathcal{G}^{\delta}$-optimal
strategy $\pi_{g^{\delta}\left(  \mathbf{m}+\mathbf{1}\right)  }^{\delta}$.
Hence, we obtain the second inequality from
\[
v^{\delta}\left(  \mathbf{m}+\mathbf{1}\right)  e^{-(c+\lambda)\delta}\leq
T_{0}(v^{\delta})\left(  \mathbf{m}\right)  \leq T(v^{\delta})\left(
\mathbf{m}\right)  =v^{\delta}\left(  \mathbf{m}\right)  .
\]

(2) In order to avoid any confusion, in the remainder of the proof we put a
superindex $\delta$ to the control actions in $\mathcal{G}^{\delta}$\textbf{.
}Note first that given any surplus in $\mathbf{R}_{+}^{n}$, the strategy of
paying dividends in such a way that the surplus goes to the nearest smaller
point in $\mathcal{G}^{2\delta}$ corresponds to go first to the nearest
smaller point in $\mathcal{G}^{\delta}$ and then to apply (possibly) a
combination of control actions $\mathbf{E}_{i}^{\delta\prime}s $. Consider
$\pi_{2g^{\delta}(\mathbf{m})}\in\Pi_{g^{2\delta}\left(  \mathbf{m}\right)
}^{2\delta}$ given by the random sequence $\mathbf{s}=(s_{k})_{k=1,...,\tilde
{k}}$ with
\[
s_{k}\in\mathcal{E}^{2\delta}=\{\mathbf{E}_{s}^{2\delta},\left(
\mathbf{E}_{i}^{2\delta}\right)  _{i=1,...,n},\mathbf{E}_{0}^{2\delta}\}.
\]
We can see that $\pi_{2g^{\delta}(\mathbf{m})}$ also belongs to $\Pi
_{2g^{\delta}(\mathbf{m})}^{\delta}$ rewriting the sequence as follows: If
$s_{k}=\mathbf{E}_{i}^{2\delta}$, we replace it by the pair $\mathbf{E}%
_{i}^{\delta}\mathbf{,E}_{i}^{\delta}$; if $s_{k}=\mathbf{E}_{s}^{2\delta}
$\textbf{,} we replace it by $\mathbf{E}_{s}^{\delta}$; and if $s_{k}%
=\mathbf{E}_{0}^{2\delta}$, we replaces it

\begin{itemize}
\item either by $\mathbf{E}_{0}^{\delta},\mathbf{E}_{0}^{\delta}$ if the next
jump in the uncontrolled process arrives at time $\tau>$ $2\delta;$

\item or by $\mathbf{E}_{0}^{\delta},\mathbf{E}_{0}^{\delta}\mathbf{,}$ and a
possible combination of $\mathbf{E}_{i}^{\delta\prime}s$, if it arrives at
time $\tau\in$ $(\delta$,$2\delta]$, so the surplus goes to the nearest
smaller point in $\mathcal{G}^{2\delta}$;

\item or by $\mathbf{E}_{0}^{\delta}$\textbf{, }and a possible combination of
$\mathbf{E}_{i}^{\delta\prime}s$, if it arrives at time $\tau\leq$ $\delta$,
so again the surplus goes to the nearest smaller point in $\mathcal{G}%
^{2\delta}$.

So we have the result. $\blacksquare$
\end{itemize}

\textit{Proof of Lemma \ref{Lipschitz Inequality u barra}}. Let us first prove
that
\begin{equation}%
\begin{array}
[c]{l}%
\left\vert V^{\delta_{k}}(\mathbf{y})-V^{\delta_{k}}(\mathbf{x})\right\vert \\
\leq\frac{2}{\hat{p}}V^{\delta_{k}}(\left\langle \mathbf{x}\vee\mathbf{y}%
\right\rangle ^{\delta_{k}})(\frac{e^{(c+\lambda)\delta_{k}}-1}{\delta_{k}%
})\left\Vert \left\langle \mathbf{y}\right\rangle ^{\delta_{k}}-\left\langle
\mathbf{x}\right\rangle ^{\delta_{k}}\right\Vert _{1}+2\delta_{k}%
\mathbf{a}\cdot\mathbf{p},
\end{array}
\label{Lipschitz V delta}%
\end{equation}
for any $\mathbf{x}$ and $\mathbf{y}$ in $\mathbf{R}_{+}^{n}$. Let us assume
first that $\mathbf{y}>\mathbf{x}$. We have from Lemma
\ref{v_delta es lipschitz y monotona en el reticulado},%
\[
V^{\delta_{k}}(g^{\delta_{k}}\left(  \mathbf{m}+\mathbf{e}_{i}\right)
)-V^{\delta_{k}}(g^{\delta_{k}}(\mathbf{m}))\leq V^{\delta_{k}}(g^{\delta_{k}%
}\left(  \mathbf{m}+\mathbf{1}\right)  )-V^{\delta_{k}}(g^{\delta_{k}%
}(\mathbf{m}))\leq V^{\delta_{k}}(g^{\delta_{k}}(\mathbf{m}))(e^{(c+\lambda
)\delta_{k}}-1).
\]
Let us call $\mathbf{m}_{\mathbf{y}}=\rho^{\delta_{k}}(\mathbf{y})$ and
$\mathbf{m}_{\mathbf{x}}=\rho^{\delta_{k}}(\mathbf{x})$. Then,%

\[%
\begin{array}
[c]{lll}%
V^{\delta_{k}}(\mathbf{y})-V^{\delta_{k}}(\mathbf{x}) & \leq & V^{\delta_{k}%
}(g^{\delta_{k}}(\mathbf{m}_{\mathbf{y}}))-V^{\delta_{k}}(g^{\delta_{k}%
}\left(  \mathbf{m}_{\mathbf{x}}\right)  )+\mathbf{a}\cdot(\mathbf{y}%
-g^{\delta_{k}}(\mathbf{m}_{\mathbf{y}}))\\
& \leq & (\frac{e^{(c+\lambda)\delta_{k}}-1}{\delta_{k}})V^{\delta_{k}%
}(\mathbf{y})\sum_{i=1}^{n}\frac{g_{i}^{\delta_{k}}\left(  \mathbf{m}%
_{\mathbf{y}}-\mathbf{m}_{\mathbf{x}}\right)  }{p_{i}}+\delta_{k}%
\mathbf{a}\cdot\mathbf{p}\\
& \leq & \left(  \frac{e^{(c+\lambda)\delta_{k}}-1}{\hat{p}\delta_{k}}\right)
V^{\delta_{k}}(\mathbf{y})\left\Vert g^{\delta_{k}}\left(  \mathbf{m}%
_{\mathbf{y}}-\mathbf{m}_{\mathbf{x}}\right)  \right\Vert _{1}+\delta
_{k}\mathbf{a}\cdot\mathbf{p}.
\end{array}
\]
Let us consider now $\mathbf{x}$ and $\mathbf{y}$ in $\mathbf{R}_{+}^{n}$,
consider \ $\mathbf{m}_{0}=\rho^{\delta_{k}}(\mathbf{x}\wedge\mathbf{y})$,
\[%
\begin{array}
[c]{l}%
\left\vert V^{\delta_{k}}(\mathbf{y})-V^{\delta_{k}}(\mathbf{x})\right\vert \\%
\begin{array}
[c]{ll}%
\leq & V^{\delta_{k}}(\mathbf{y})-V^{\delta_{k}}(\mathbf{x}\wedge
\mathbf{y})+V^{\delta_{k}}(\mathbf{x})-V^{\delta_{k}}(\mathbf{x}%
\wedge\mathbf{y})\\
\leq & \frac{1}{\hat{p}}V^{\delta_{k}}(\mathbf{x}\vee\mathbf{y})(\frac
{e^{(c+\lambda)\delta_{k}}-1}{\delta_{k}})\left(  \left\Vert g^{\delta_{k}%
}\left(  \mathbf{m}_{\mathbf{y}}-\mathbf{m}_{0}\right)  \right\Vert
_{1}+\left\Vert g^{\delta_{k}}\left(  \mathbf{m}_{\mathbf{x}}-\mathbf{m}%
_{0}\right)  \right\Vert _{1}\right)  +2\delta_{k}\mathbf{a}\cdot\mathbf{p}\\
\leq & \frac{2}{\hat{p}}V^{\delta_{k}}(\mathbf{x}\vee\mathbf{y})(\frac
{e^{(c+\lambda)\delta_{k}}-1}{\delta_{k}})\left\Vert g^{\delta_{k}}\left(
\mathbf{m}_{\mathbf{y}}-\mathbf{m}_{\mathbf{x}}\right)  \right\Vert
_{1}+2\delta_{k}\mathbf{a}\cdot p.
\end{array}
\end{array}
\]
Therefore we have (\ref{Lipschitz V delta}).

By definitions (\ref{Definicion Vdelta}) and (\ref{ubarra como limite}), and
since $T_{i}\left(  v^{\delta_{k}}\right)  \leq v^{\delta_{k}}$,
\[%
\begin{array}
[c]{lll}%
\overline{V}(\mathbf{y})-\overline{V}(\mathbf{x}) & \geq & \overline
{V}(\mathbf{y})-V^{\delta_{k}}(\mathbf{y})+\mathbf{a}\cdot g^{\delta_{k}%
}\left(  \rho^{\delta_{k}}(\mathbf{y})-\rho^{\delta_{k}}(\mathbf{x})\right) \\
&  & +\mathbf{a}\cdot(\mathbf{y}-g^{\delta_{k}}(\rho^{\delta_{k}}%
(\mathbf{y})-\rho^{\delta_{k}}(\mathbf{x}))+\mathbf{x})+V^{\delta_{k}%
}(\mathbf{x})-\overline{V}(\mathbf{x});
\end{array}
\]
taking the limit as $k$ goes to infinity, we obtain the first inequality of
the Lipschitz inequality.

We can write, from (\ref{Lipschitz V delta}),
\[%
\begin{array}
[c]{lll}%
\overline{V}(\mathbf{y})-\overline{V}(\mathbf{x}) & = & \overline
{V}(\mathbf{y})-V^{\delta_{k}}(\mathbf{y})+V^{\delta_{k}}(\mathbf{y}%
)-V^{\delta_{k}}(\mathbf{x})+V^{\delta_{k}}(\mathbf{x})-\overline
{V}(\mathbf{x})\\
& \leq & \overline{V}(\mathbf{y})-V^{\delta_{k}}(\mathbf{y})+\frac{2}{\hat{p}%
}\overline{V}(\mathbf{y})(\frac{e^{(c+\lambda)\delta_{k}}-1}{\delta_{k}%
})\left\Vert g^{\delta_{k}}\left(  \rho^{\delta_{k}}(\mathbf{y})-\rho
^{\delta_{k}}(\mathbf{x})\right)  \right\Vert _{1}\\
&  & +2\delta_{k}\mathbf{a}\cdot\mathbf{p}+V^{\delta_{k}}(\mathbf{x}%
)-\overline{V}(\mathbf{x});
\end{array}
\]
taking the limit as $k$ goes to infinity, we obtain the second inequality of
the Lipschitz inequality.$~\blacksquare$


\begin{thebibliography}{99}                                                                                               %
\bibitem {AlAZMU}Albrecher, H., Azcue, P., and Muler, N. (2017). Optimal
dividend strategies for two collaborating insurance companies. \textit{Adv.
App. Prob. } 49(2), 515-548.

\bibitem {Asmussen Taksar 1997}Asmussen S. and Taksar M. (1997). Controlled
diffusion models for optimal dividend pay-out. \textit{Insurance Math.
Econom.} 20, 1--15.

\bibitem {APP 2007}Avram, F., Palmowski Z. and Pistorius, M. (2007). On the
optimal dividend problem for a spectrally negative L\'{e}vy process.
\textit{Ann. Appl. Probab.} 17, 156-180.

\bibitem {APP 2015}Avram, F., Palmowski Z. and Pistorius, M. (2015). On
Gerber-Shiu functions and optimal dividend distribution for a Levy
risk-process in the presence of a penalty function. \textit{Ann. Appl.
Probab.} 25 (4), 1868--1935

\bibitem {AM 2005}Azcue, P. and Muler N.\textsc{\ }(2005). Optimal reinsurance
and dividend distribution policies in the Cram\'{e}r-Lundberg model\textit{.
Math. Finance }15, 261-308.

\bibitem {AM Libro}Azcue P. and Muler N. (2014). \textit{Stochastic
Optimization in Insurance: a Dynamic Programming Approach. }Springer Briefs in
Quantitative Finance. Springer.

\bibitem {AM Switching}Azcue P. and Muler N. (2015). Optimal dividend payment
and regime switching in a Compound Poisson risk model. \textit{SIAM J. Control
Optim.} 53(5), 3270-3298.

\bibitem {BS}Barles G. and Souganidis P. E. (1991). Convergence of
approximation schemesfor fully nonlinear second order equations.
\textit{Asymptot. Anal.} 4:271--283.

\bibitem {BR}Budhiraja A. and Ross K. (2007). Convergent numerical scheme for
singular stochastic control with state constraints in a portfolio selection
problem. \textit{SIAM J. Control Optim.}, 45(6):2169--2206.

\bibitem {CL}Crandall, M. G. and Lions, P. L.\textsc{\ }(1983). Viscosity
solution of Hamilton-Jacobi equations. \textit{Trans. Amer. Math. Soc.} 277, 1-42.

\bibitem {CP}Czarna I. and Palmowski, Z. (2011). De Finetti's dividend problem
and impulse control for a two-dimensional insurance risk process.
\textit{Stochastic Models} \textbf{27}, 220--250.

\bibitem {De Finetti}De Finetti, B. (1957). Su Un'impostazione Alternativa
Della Teoria Collettiva Del Rischio. \textit{Transactions of the XVth
International Congress of Actuaries} 2 , 433--443.

\bibitem {DicksonWaters2004}Dickson, D.C.M. and Waters, H.R. (2004). Some
optimal dividends problems. \textit{ASTIN Bulletin} 34, 49--74.

\bibitem {FS}Fleming, W.H. and Soner, H.M. (1993). \textit{Controlled Markov
Processes and Viscosity Solutions}. New York: Springer-Verlag.

\bibitem {Gerber}Gerber, H. (1969). Entscheidungskriterien f\"{u}r den
zusammengesetzten Poisson-Proze\ss , \textit{Mitt. Ver. Schweiz. Vers. Math}
69, 185-228.

\bibitem {GerberLinYang2006}Gerber H. U., Lin X. S. and Yang H (2006). A Note
on the Dividends-Penalty Identity and the Optimal Dividend Barrier.
\textit{Astin Bulletin }36, 489-503.

\bibitem {GS 1998}Gerber H.U., Shiu E.S.W (1998). On the time value of ruin.
\textit{N. Am. Actuar. J.} 2 (1), pp. 48-78.

\bibitem {GS Merger}Gerber H. U. and Shiu E. S. W. (2006). On the merger of
two companies. \textit{N. Am. Actuar. J.} 10(3): 60-67.

\bibitem {KD}Kushner H. J., and Dupuis P. (2001). \textit{Numerical Methods
for Stochastic Control Problems in Continuous Time}. Springer Science \&
Business Media.

\bibitem {KM}Kushner H.J. and Martins L.F. (1991). Numerical methods for
stochastic singular control problems. \textit{SIAM J. Control Optim.} 29:1443--1475.

\bibitem {Loeffen2008}Loeffen, R.L. (2008). On optimality of the barrier
strategy in de Finetti's dividend problem for spectrally negative Levy
processes. \textit{Ann. Appl. Probab.} 18, 1669--1680.

\bibitem {LoeffenRenaud2010}Loeffen, R. L. and Renaud, J. F. (2010). De
Finetti's optimal dividends problem with an affine penalty function at ruin.
\textit{Insurance Math. Econom.} 46, 98-108.

\bibitem {LPV}Ly Vath, V., Pham, H. and Villeneuve S. (2008). A mixed
singular/switching control problem for a dividend policy with reversible
technology investment. \textit{Ann. Appl. Probab} 18, 1164--1200.

\bibitem {Schmidli book 2008}Schmidli H. (2008). \textit{Stochastic Control in
Insurance}. Springer, New York.

\bibitem {S}Souganidis P. E. \ (1985). Approximation schemes for viscosity
solutions of Hamilton-Jacobi equations. \textit{Jour. Differential Eqns.} 59, 1-43.

\bibitem {TA}Thonhauser, S. and Albrecher, H. (2007). Dividend maximization
under consideration of the time value of ruin. \textit{Insurance Math.
Econom.} 41 163--184.
\end{thebibliography}
\end{document}